\documentclass[12pt]{article}
\usepackage[margin=1in]{geometry}
\usepackage{amsmath, amsfonts, amssymb, amsthm, mathrsfs, dsfont}
\usepackage{graphicx}
\usepackage[font=small,labelfont=bf]{caption}
\usepackage{tikz}
\usetikzlibrary{arrows.meta,calc,decorations.markings,math,arrows.meta}
\tikzset{dot/.style={fill=black,circle}}

\begin{document}
\title{Solving Maxwell's Equation in 2D with NNLCI}
\author{Harris Cobb, Hwi Lee, Yingjie Liu}
\date{}
\maketitle

\newcommand{\lb}{\langle} \newcommand{\rb}{\rangle}
\newcommand{\ls}{\lesssim} \newcommand{\gs}{\gtrsim}
\newcommand{\nls}{\lnsim} \newcommand{\ngs}{\gnsim}

\begin{abstract}
In this paper we apply neural networks with local converging inputs (NNLCI) , originally developed in \cite{hly21}, to solve the two dimensional Maxwell's equation around perfect electric conductors (PECs). The input to the networks consist of local patches of low cost numerical solutions to the equation computed on two coarse grids, and the output is a more accurate solution at the center of the local patch. We apply the recently developed second order finite difference method \cite{lee2022ghost} to generate the input and training data which captures the scattering of electromagnetic waves off of a PEC at a given terminal time. The advantage of NNLCI is that once trained it offers an efficient alternative to costly high-resolution conventional numerical methods; our numerical experiments indicate the computational complexity saving by a factor of $8^3$ in terms of the number of spatial-temporal grid points. In contrast with existing research work on applying neural networks to directly solve PDEs, our method takes advantage of the local domain of dependence of the Maxwell's equation in the input solution patches, and is therefore simpler, yet still robust. We demonstrate that we can train our neural network on some PECs to predict accurate solutions to different PECs with quite different geometries from any of the training examples.
\end{abstract}

\begin{section}{Introduction}
In recent years, applications of artificial neural networks have permeated the broad realm of computational science thanks to the advances in computing powers such as Graphical Processing Units (GPUs) as well as the expansion of data availability. Some specific examples encompass biological applications \cite{ndg20}, image processing \cite{ksh17}, uncertainty quantification \cite{tb18} among a plethora of existing work. The utility of the neural networks (NNs) can be attributed to the classical result on its universal approximation property \cite{c89, hsw90}, of which our current work is yet another demonstration in the context of numerically solving partial differential equations (PDEs). More specifically we are interested in the so-called forward problem of computing numerical solutions as opposed to the backward problems of identifying the underlying PDEs, although the latter has received as much attention in scientific computing literature  \cite{meidani2021data,kang2021ident,he2020robust}.

An approach to using NNs to solve PDEs is to approximate the mapping from all known information (such as the initial value, boundary values etc) to the unknown solution by a deep neural network and use many existing solutions for training, see e.g \cite{lkal20}. Another approach is to approximate the mapping from a space-time location to the solution there with a deep neural network and define the loss function as certain form of the residual error of a PDE, plus penalty terms enforcing the physical data such as initial and boundary conditions. The deep Galerkin method \cite{sirignano2018dgm} introduces a loss function which contains the residue error of a finite difference approximation to the PDE. The deep Ritz method \cite{EYu18} introduces a loss function which contains the Ritz energy of the finite element approximation of the PDE. The physics informed neural network (PINN) \cite{raissi2019physics} exploits the exact auto-differentiation \cite{baydin2018automatic} in computation of the loss function which is designed to enforce consistency with the governing PDE. These methods harness the approximation property of NNs as an ansatz for the solutions to the PDEs. Besides the solution learning, researchers have also investigated the possibilities of training NNs to learn the solution operators as a function of the parameters of the PDE \cite{zzpp19}.
The domains for successful application of PINN are far-reaching, including but not limited to fluid flows, heat transfer and optics \cite{mjk20,cwwpk21,clkd20}. The ongoing research efforts have resulted in some modifications of PINN itself as in \cite{plk19}, incorporating information from partially known solutions in low and high fidelities in the loss function \cite{meng2020composite} and  the emergence of related approaches like physics-informed convolutional neural networks used in \cite{zgzyc23}. We would like to refer to \cite{cds22} and the references therein for an overview of PINN and more classical NN-based approaches. We also point out some alternatives such as those based on Gaussian processes \cite{pang2019neural}, and NN augmented time-stepping strategies \cite{siahkoohi1910neural}.

The accuracy enhancement by NNs has been illustrated in existing literature  \cite{pcg12}, and the growing number of studies is geared towards complexity-reduction hybrid methods combining NNs and conventional numerical methods  \cite{liu2019multi, nguyen2021numerical}.
In this work we pursue a different NN-based approach, adopting the architecture of recently developed neural networks with local converging inputs (NNLCIs) \cite{hly21,hly22}.  The NNLCIs have been successfully applied to nonlinear Euler equations for gas dynamics in one and two spatial dimensions as a black-box tool to generate high fidelity solutions containing complicated interactions of shocks and contacts from two converging, low-cost, low-fidelity numerical solutions.  Our work pursues another application of the NNLCIs as a computational saving tool, instead of investigating the high expressivity and intricate inner workings of which have yet to be fully and rigorously understood. The onus is on the NNs to acquire the predictive ability for accuracy improvement when the users work on careful design of the inputs to the network inputs rather than the network architectures such as choices for activation functions.

We demonstrate the utility of NNLCIs in the case of solving Maxwell's equations around curved perfect electric conductors with corners. The scattering of electromagnetic waves is of practical importance with industrial applications as in antenna designs \cite{mlcv12}, for which one may need to do repeated computations with various shapes of PEC objects. It is indeed the high complexity of possible variations in PEC boundary geometries that distinguishes our work from the earlier applications of NNLCIs. We illustrate that it is sufficient for the NNLCIs to be trained locally with respect to some PEC shapes to be able to extract more accurate solutions from the two coarse grid solutions for different shapes of PECs. The parameter space of local geometries of PECS is expected to be more tractable than that of global ones, and after all the analytical, hence numerical, solutions to the Maxwell's equations have local domains of dependence. We fully exploit such locality in approximation of the underlying unknown multivariate mapping, steering away from more conventional approach of taking as a single input all the known data from the entire computational domain.

Our proposed NNLCI-based approach benefits largely from the recently developed numerical scheme \cite{hly22} which is used to compute the inputs to our network. What makes the scheme attractive is its simple and systematic treatment of boundary conditions for general curved PEC shapes that may even have corners. More specifically the scheme is based on automatic construction of ghost values by means of the level set framework \cite{OsherSethian88}, the PDE-based extension technique \cite{Peng99,Fedkiw99} and the so-called guest values. Another merit of the scheme is its second order accuracy which is what the popular Yee scheme \cite{y96} can achieve in free space without any PECs. 
We also point out that our underlying numerical scheme is based on uniform rectangular grids, which are particularly well-suited to NNLCIs since the inputs are finite difference stencils at two different resolutions that need to be aligned at the same spatial and temporal locations.

The rest of the paper is organized as follows. In Section \ref{sec:2}, we introduce the scattering problem and the numerical scheme we use to produce the training and testing data. We then present our NNLCI focusing on how its inputs are to be formatted. In Section \ref{sec:3}, we provide the results of our numerical experiments to illustrate the effectiveness of our approach. We conclude our paper in Section \ref{sec:4}.

\end{section}

\begin{section}{Problem Formulation}
\label{sec:2}
\begin{subsection}{The governing PDEs and underlying numerical method}
The classical Maxwell's equations govern the dynamics of electromagnetic fields through a medium which we assume in this work is isotropic and homogeneous. In particular we consider the scattering of electromagnetic waves off of a PEC in free space, that is, 
$$
\begin{aligned}
\epsilon_r\frac{\partial \bf{E}}{\partial t} &=\nabla \times \bf{H} \\
\mu_r\frac{\partial \bf{H}}{\partial t} &=-\nabla \times \bf{E}
\end{aligned}
$$
subject to some initial conditions and the PEC boundary conditions that are given by
$$\bf{E}\times \bf{n}=\bf{0} \text{ and } \bf{H} \cdot \bf{n}=0.$$
Here $\epsilon_r$ and $\mu_r$ denote the relative permittivity and permeability, respectively, and $\bf{n}$ denotes the outward unit normal vector to the PEC interface. Without loss of generality we set $\epsilon_r = \mu_r = 1$, and assume for simplicity that the waves are transmagnetic polarized in the $z-$direction, thereby obtaining the following two dimensional system of equation
\begin{eqnarray}
\label{eqn:mw}
\frac{\partial H_x}{\partial t} &=&-\frac{\partial E_z}{\partial y}  \nonumber \\
\frac{\partial H_y}{\partial t} &=&\frac{\partial E_z}{\partial x}  \\
\frac{\partial E_z}{\partial t}&=&\frac{\partial H_y}{\partial x}-\frac{\partial H_x}{\partial y}. \nonumber
\end{eqnarray}
In order to numerically simulate the dynamics (\ref{eqn:mw}) in free space, we first choose the unit square $[0,1]^2$ as our computational domain. We superimpose on the domain  uniform rectangular grids $\{(x_i,y_j)\}$ with $\Delta x = \Delta y$. We adopt the implementation of  perfectly matched layer as done in \cite{s10} to enforce the far field absorbing boundary conditions so that the reflected waves off the PEC do not re-enter the computational domain. To this end we follow the total/scatter formulation \cite{im11}  with the prescription of the incident waves.

We apply the finite difference method \cite{lee2022ghost} to compute the numerical solutions to (\ref{eqn:nsmw}), which will serve as training and test data to our NNLCI. The method is mainly composed of three elements, the first of which is the first order scheme 
\begin{equation}
\label{eqn:nsmw}
U^{n+1}_{i,j} = \left(\frac{1}{5}\mathcal{C}+\frac{4}{5}\mathcal{F}\right)U^{n}_{i,j}.
\end{equation} where $U^n_{i,j}$ denotes the numerical scheme at the time level $t = n\Delta t $ and at the spatial location $(x_i,y_j)$. Here $\mathcal{C}$ denotes the forward-in time, centered-in-space numerical scheme while $\mathcal{F}$ denotes the classical Lax-Friedrichs scheme. We choose scheme (\ref{eqn:nsmw}) as it has already been used in \cite{lee2022ghost,zl22} to simulate the scattering of electromagnetic waves about PECs, but one may instead consider a more general scheme $\left(\theta \mathcal{C}+(1-\theta) \mathcal{F}\right)$ with $0\leq \theta \leq 1$ \cite{w19}. The second element is the construction of ghost values inside the PECs via a level-set based PDE extension technique as well as the introduction of guest values. The constructed ghost values are to be locally second order accurate to ensure the first order accuracy of the scheme \eqref{eqn:nsmw}. The last ingredient is the Back and Forth Error Compensation and Correction (BFECC) method \cite{dupont2003back} which in general increases the order of accuracy by one for odd order schemes \cite{owhadi2019kernel}. Therefore in our current setting our numerical solutions are computed to second order accuracy. Moreover BFECC provides enhanced numerical stability by allowing a larger CFL number, hence we choose  $\Delta t = \Delta x$ which would otherwise not be allowed.  An example of the computed numerical solutions is shown in Figure \ref{fig:snapshot} which contains only the total field region.
\begin{center}
\includegraphics[scale=0.5]{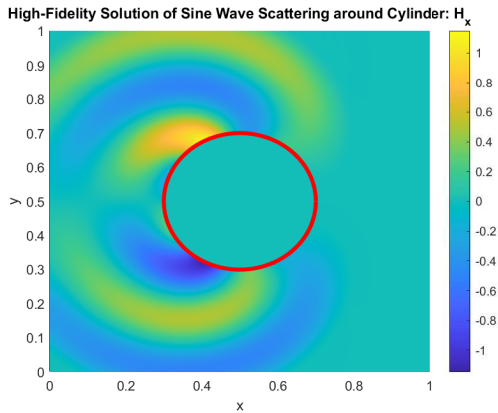}  
\captionof{figure}{$H_x$ scattering around a cylinder after 0.8 seconds.} \label{fig:snapshot}
\end{center}
 
The question of interest to us is whether  it is possible for a neural network to produce high fidelity numerical solutions from numerical solutions computed on two coarse grids. The affirmative answer, at least in the cases we demonstrate, is made possible by the NNLCI which we present in the sequel.

\end{subsection}

\begin{subsection}{Application of NNLCIs}

We follow the architecture of the neural network used in \cite{hly21,hly22} which is a simple fully connected feed-forward network based on multiple steps of Adams optimizer and a final step of L-BFGS for each iteration instead of stochastic gradient descent. We modify the hyperparameters of the optimizer used in their work so that the learning rate of L-BFGS is 3e-6 and the learning rate of Adams is 1e-5. 
Our network consists of 4 layers with 50 neurons in each layer, hence it can be trained readily yet we observe highly predictive capabilities of the network in our numerical experiments.

We iterate the numerical scheme introduced above until a fixed terminal time, $T = 0.8$ in most of our experiments unless stated otherwise, on three different numerical resolutions, two low-fidelity resolutions and one high-fidelity resolution. In particular, we take $\Delta x = 1/40, 1/80$ for the two low fidelity solutions and $\Delta x = 1/640$ for the high-fidelity solution which we regard as the ``exact" reference solution. 
We store all the numerical solutions only at the terminal time, as opposed to the two different time levels as done in the earlier work \cite{hly21}. This is because the Maxwell's equations (\ref{eqn:mw}) are time reversible, hence the solutions at the previous time step are uniquely determined by the solutions at the current time step, so we omit any notation that denotes the time. 

From the stored numerical solutions, we prepare the input data to our network as follows. Let us write $(i',j'), (i'',j''),$ and $(i,j)$ to denote the same location in space for the coarsest, finer and the reference grids respectively. We consider a $3\times 3$ window in the coarsest grid centered at $(i',j')$. In this paper, we restrict to the case where the center, $(i',j')$, and all its 8 neighboring grid points belong to the total-field region outside the PEC. In the future, we will consider relaxation of this restriction. 
From the coarsest grid we reshape the $3\times 3$ windows for each of $H_x,H_y,E_z$ into $1\times 9$ vectors 
$$\{(H_x)_{i'-1,j'-1},(H_x)_{i',j'-1},(H_x)_{i'+1,j'-1},\ldots,(H_x)_{i'-1,j'+1},(H_x)_{i',j'+1},(H_x)_{i'+1,j'+1}\}$$
$$\{(H_y)_{i'-1,j'-1},\ldots,(H_y)_{i'+1,j'+1}\}, \{(E_z)_{i'-1,j'-1},\ldots,(E_z)_{i'+1,j'+1}\}$$ and concatenate all $27$ of the values into a single vector. Similarly we take from the finer grid $$\{(H_x)_{i''-2,j''-2},(H_x)_{i'',j''-2},(H_x)_{i''+2,j''-2},\ldots,(H_x)_{i''-2,j''+2},(H_x)_{i'',j''+2},(H_x)_{i''+2,j''+2}\}$$ $$\{(H_y)_{i''-2,j''-2},\ldots,(H_y)_{i''+2,j''+2}\}, \{(E_z)_{i''-2,j''-2},\ldots,(E_z)_{i''+2,j''+2}\}$$ and concatenate all the values. As a result the input vector consists of 54 input values $$\{(H_x)_{i'-1,j'-1},\ldots,(E_z)_{i''+2,j''+2}\},$$ which is then to be fed into our network to produce the output vector $$\{(H_x)_{i,j},(H_y)_{i,j},(E_z)_{i,j}\}.$$ We illustrate how the input data are prepared in Fig.~\ref{fig:input} where $\bf{U}$ represents the vector containing all three components $H_x,H_y,E_z$.

\begin{tikzpicture}
\draw[step=2cm,black] (0,0) grid (4,4);
\draw (0,4.5) node{$\bf{U}_{i'-1,j'+1}$};
\draw (2,4.5) node{$\bf{U}_{i',j'+1}$};
\draw (4,4.5) node{$\bf{U}_{i'+1,j'+1}$};
\draw (-0.75,2.5) node{$\bf{U}_{i'-1,j'}$};
\draw (1.25,2.5) node{$\bf{U}_{i',j'}$};
\draw (3.25,2.5) node{$\bf{U}_{i'+1,j'}$};
\draw (-0.75,0.5) node{$\bf{U}_{i'-1,j'-1}$};
\draw (1.25,0.5) node{$\bf{U}_{i',j'-1}$};
\draw (3.25,0.5) node{$\bf{U}_{i'+1,j'-1}$};
\draw [-Latex, line width = 1.5mm](5,2) -- (10,0.1);
\draw (8.5,3) node{Low-Fidelity Input 1};
\draw [-Latex, line width = 0.75mm](6.25,3) -- (5,2.5);
\draw (10.5,-0.5) node{$\{(H_x)_{i'-1,j'-1},(H_x)_{i',j'-1},\ldots,(E_z)_{i''+2,j''+2}\}$};
\node[dot] at (0,4){};
\node[dot] at (2,4){};
\node[dot] at (4,4){};
\node[dot] at (0,2){};
\node[dot] at (2,2){};
\node[dot] at (4,2){};
\node[dot] at (0,0){};
 \node[dot] at (2,0){};
\node[dot] at (4,0){};

\draw (8.5,-4) node{Low-Fidelity Input 2};
\draw [-Latex, line width = 0.75mm](6.25,-4) -- (5,-3.5);
\draw[step=1cm,black] (0,-5) grid (4,-1);
\draw (0,-0.5) node{$\bf{U}_{i''-2,j''+2}$};
\draw (2,-0.5) node{$\bf{U}_{i'',j''+2}$};
\draw (4,-0.5) node{$\bf{U}_{i''+2,j''+2}$};
\draw (-0.75,-2.5) node{$\bf{U}_{i''-2,j''}$};
\draw (1.25,-2.5) node{$\bf{U}_{i'',j''}$};
\draw (3.25,-2.5) node{$\bf{U}_{i''+2,j''}$};
\draw (-0.75,-4.5) node{$\bf{U}_{i''-2,j'-2}$};
\draw (1.25,-4.5) node{$\bf{U}_{i'',j''-2}$};
\draw (3.25,-4.5) node{$\bf{U}_{i''+2,j''-2}$};
\draw [-Latex, line width = 1.5mm](5,-3) -- (10,-1.1);
\node[dot] at (0,-1){};
\node[dot] at (2,-1){};
\node[dot] at (4,-1){};
\node[dot] at (0,-3){};
\node[dot] at (2,-3){};
\node[dot] at (4,-3){};
\node[dot] at (0,-5){};
\node[dot] at (2,-5){};
\node[dot] at (4,-5){};

\end{tikzpicture}

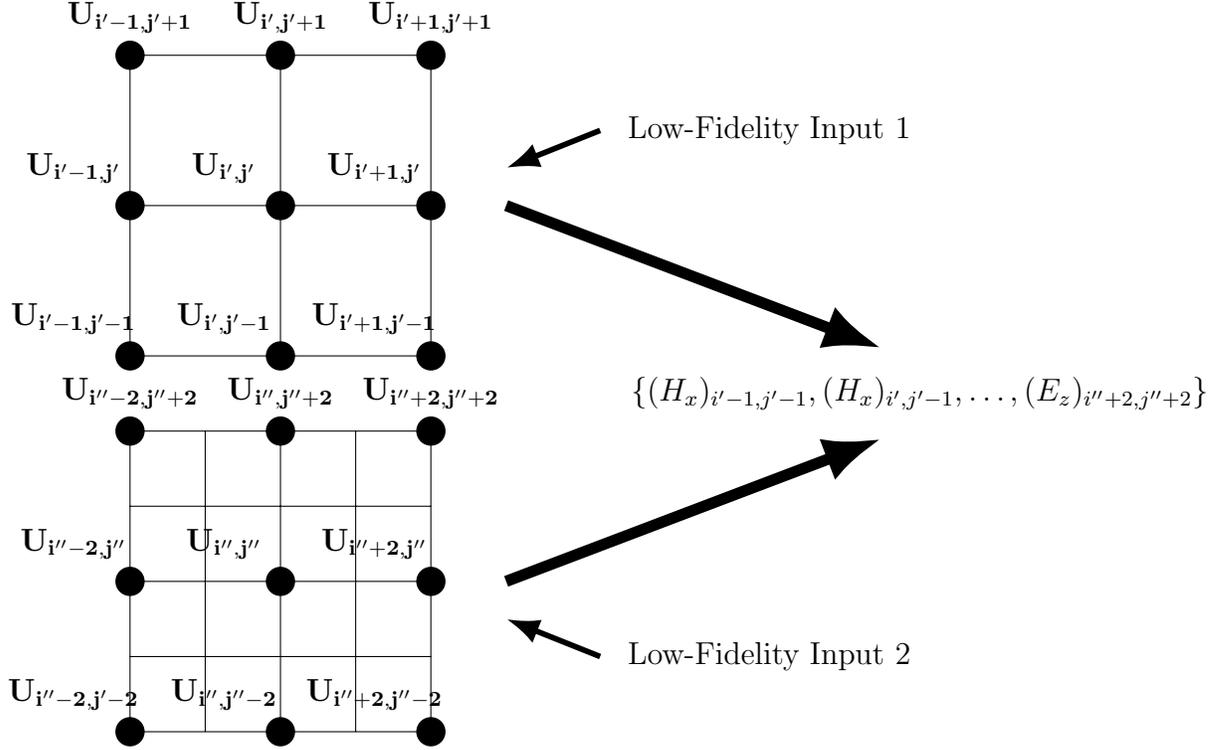
\captionof{figure}{Procedure for preparing input to 2CGNN.}\label{fig:input}

The ordering of the input data is immaterial so far as it is consistent across all $3 \times 3$ windows.
Following the approach in \cite{hly22}, one may include the time step size $\Delta t$ as part of the input data. However in the present work the simpler input format we choose is sufficient to produce accurate solutions.

We train the network until the training error is 1e-5. The training and testing data consist of perturbations of the incident wave given by
\begin{align}
\label{eqn:ww}
H_x(x,y,t) &= 0  \nonumber \\
H_y(x,y,t) &=  - \sin( \omega (x-t))) \chi(t-x)\\
E_z(x,y,t) &=  \sin( \omega (x-t)) \chi(t-x) \nonumber
\end{align} where  $\omega = 2\pi/0.3$ is the angular frequency and $\chi(z)$ denotes the heaviside step function. We also consider different shapes of the PECs to illustrate the robustness of the network.  Once the network is trained, we measure its performance by comparing the relative errors (over the total field region) of the predicted numerical solutions and the finer grid input solutions.

\end{subsection}
\end{section}

\begin{section}{Numerical Experiments}
\label{sec:3}
\begin{subsection}{Different source waves}
We first demonstrate the performance of our network to predict high-fidelity solutions when the incident waves in \eqref{eqn:ww} are perturbed with respect to their amplitudes or the angular frequencies. Unless stated otherwise, we fix the geometry of the PEC to a circular PEC centered at $(0.5,0.5)$ with the radius $0.2$. Moreover both the training and testing input vectors are taken from all $3\times 3$ admissible input windows over the entire computational domain.

\begin{subsubsection}{Varying amplitudes}
 
 We consider up to $10\%$ perturbations of the amplitude of the incident waves. More specifically we compute the two coarse grid and reference solutions subject to the perturbations of $0,\pm 2,\pm 3,\pm 4,\pm 5,\pm 6,\pm 7,\pm 8,\pm 10 \%$ corresponding to the $17$ data sets of amplitudes from $0.9$ to $1.1$. The computed data are then separated into either the training or testing sets. 
 The first experiment consists of training on the even perturbations and testing on the odd perturbations. The second is to train on $0,\pm 10 \%$  perturbations, and test {on the rest of the perturbations}. Lastly we train on $0,\pm 5 \%$ perturbations and test on the rest. In all the cases we consider throughout this paper, we decrease to the same network training loss, $2.5e-5$. We compute the relative errors in $l_2$ norm of the predicted and finer grid input solutions with respect to the reference solutions. The averaged errors with respect to the number of perturbations in the testing data are shown in Table \eqref{tbl:ap}.
\begin{center}
\begin{tabular}{ |c|c|c|c|c|c|  }
 \hline
 \multicolumn{6}{|c|}{Circular PEC: modifying amplitudes} \\
 \hline
 \multicolumn{2}{|c|}{Even Perturbations} & \multicolumn{2}{|c|}{$0\%,\pm 10\%$} & \multicolumn{2}{|c|}{$0\%,\pm 5\%$}\\
 \hline
Finer Input Err. & 0.1634 & Finer Input Err. & 0.1634 & Finer Input Err. & 0.1634\\
Pred. Err. & 0.0031 & Pred. Err. & 0.0029 & Pred. Err. & 0.0042\\
 \hline
\end{tabular} 
\captionof{table}{Accuracy of the predicted solutions to a circular PEC with varying amplitudes. }
\label{tbl:ap}
\end{center}

 Since the Maxwell's equations \eqref{eqn:mw} are scaling invariant, one may hence use a small set of training data with fairly large training gaps so far as the perturbed amplitudes are concerned. We note that in the first two experiments (training on the even perturbations, and $0\%,\pm 10\%$) all the testing vectors belong to the convex hull of the training data, hence the corresponding outputs are linearly interpolated from the training data. However the last case is concerned with some testing vectors for which the desired outputs are to be linearly extrapolated, which suggest one of the robust features of our method which was not considered in the earlier work \cite{hly21}. 
\end{subsubsection}

\begin{subsubsection}{Varying frequencies}
We next consider modifying the angular frequencies of the incident waves. It is of interest to demonstrate the performance of the network when the  effects of the perturbed frequencies are nonlinear on the corresponding solutions. We generate the training and testing data up to $ \pm 10\%$ as in the amplitude case. 
The results are displayed in Figures \ref{Fig:vfcs} and \ref{Fig:vfcc}. We observe from Table \ref{tbl:cfreq} that the accuracy gain is significant when the training gap is small, covering all the even perturbations. However, there is not much improvement in accuracy when the training gap is enlarged to $\pm 10 \%$ when the network is expected to quadratically interpolate the training data. On the other hand we do not observe better accuracy improvement when the training gap alone is narrowed to the $\pm 5\%$ perturbations without increasing the number of the training data. To rule out the possibility of limited accuracy gain due to extrapolation, we have limited the testing data only to $\pm 7 \%$ and $\pm 8 \%$ subject to two training cases, the first on $0\%,\pm 10\%$ and the second on $0\%,\pm 5\%$. We observe in Table \ref{tbl:cfreqexin} that the number of training data plays a more critical role than the training gap as it caps the complexity of the NN approximation.

Note that in following and subsequent graphs the $H_x$ field appears to be non-zero within the boundary of the PEC. This is only due to the limitation of how MatLab generates surface plots: The color of each cell is determined by the value of the top-left corner, so the interior conditions of the PEC are not violated. Also, since we only use as input $3\times 3$ windows which are entirely contained within the total field, we cannot as of yet predict values which are within $\Delta x$ of the PEC. To this end for the purposes of graphical clarity, we take the values around the PEC which cannot be predicted from the high-fidelity solution. However in future work, we will attempt to predict values all the way up to and including the boundary of the PEC which will be owed to the modification of our numerical scheme which will give a higher-order recovering of the ghost values.

\begin{center}
\includegraphics[scale=0.5]{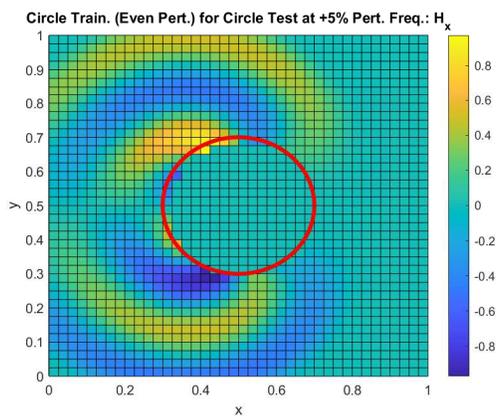}\\
\captionof{figure}{Predicted values of $H_x$ scattering off a circle for the $+5\%$ perturbed frequency in the plane at 0.8 seconds after training the network on all the even perturbations.}\label{Fig:vfcs}
\end{center}

\includegraphics[scale=0.5]{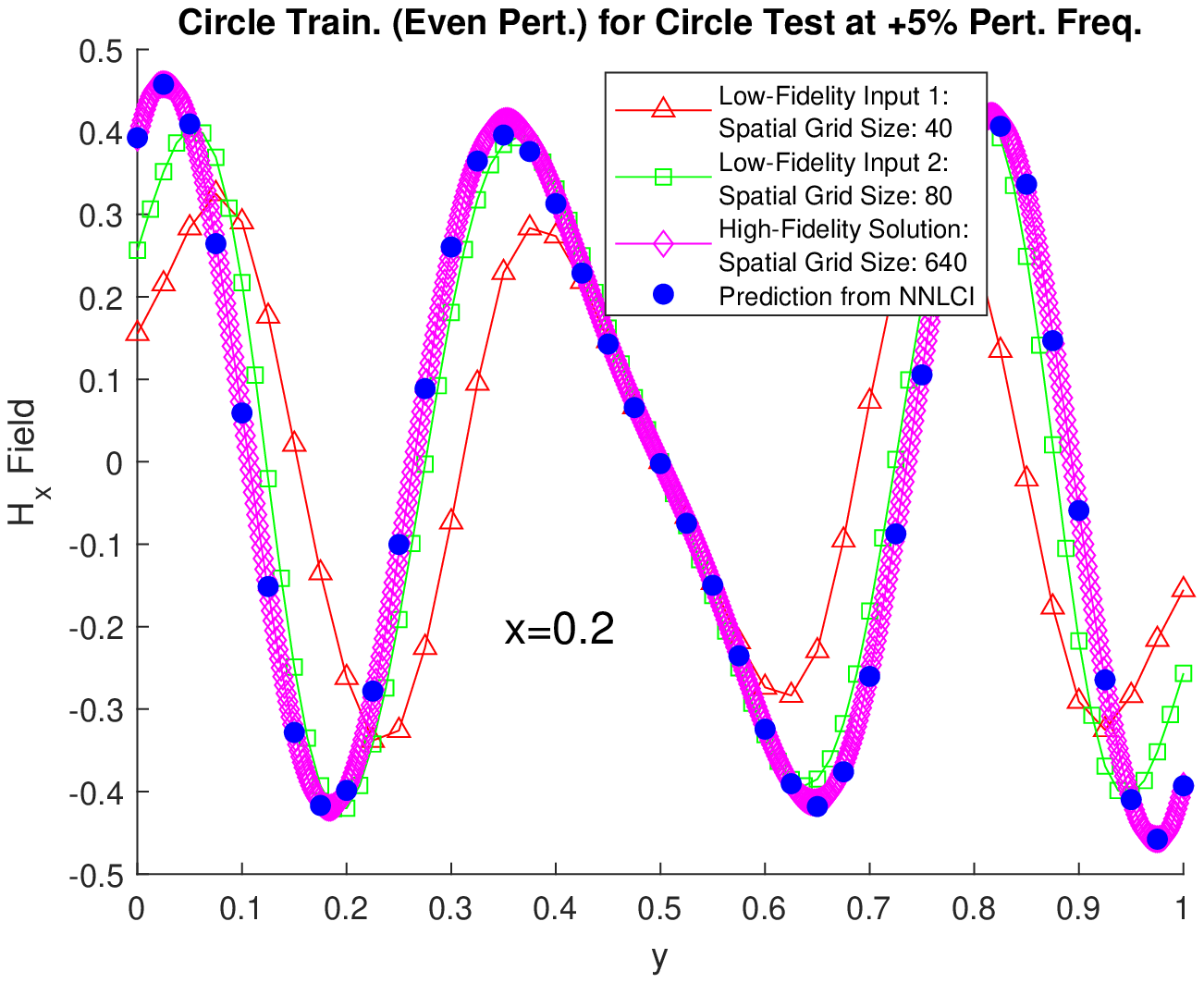}
\includegraphics[scale=0.5]{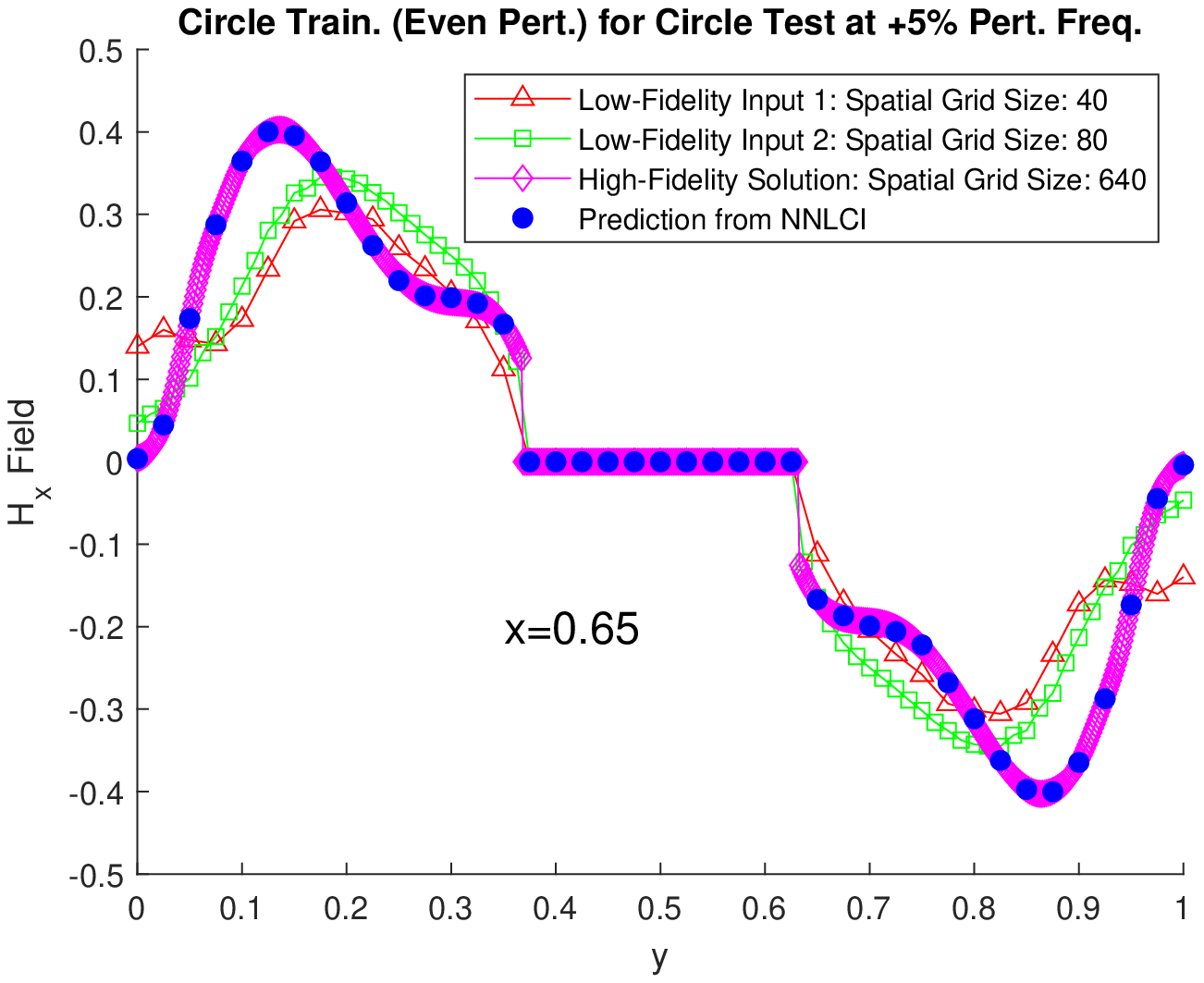}
\captionof{figure}{Cross-sections of the computed $H_x$ at different $y$ values.}\label{Fig:vfcc}

\begin{tabular}{ |c|c|c|c|c|c|  }
 \hline
 \multicolumn{6}{|c|}{Circular PEC: modifying frequency} \\
 \hline
 \multicolumn{2}{|c|}{Even Perturbations} & \multicolumn{2}{|c|}{$0\%,\pm 10\%$} & \multicolumn{2}{|c|}{$0\%,\pm 5\%$}\\
 \hline
Finer Input Err. & 0.1663 & Finer Input Err. & 0.1664 &Finer Input Err. & 0.1672\\
Pred. Err. & 0.0087 & Pred. Err. & 0.0430 & Pred. Err. & 0.0651\\
\hline
\end{tabular}
\captionof{table}{Accuracy of the predicted solutions to circular PEC with varying frequencies.}
\label{tbl:cfreq}

\begin{center}
\begin{tabular}{ |c|c|c|c|  }
 \hline
 \multicolumn{4}{|c|}{Circular PEC: modifying frequency} \\
 \hline
 \multicolumn{2}{|c|}{$0\%,\pm 10\%$} & \multicolumn{2}{|c|}{$0\%,\pm 5\%$}\\
 \hline
Finer Input Err. & 0.1690 & Finer Input Err. & 0.1690\\
Pred. Err. & 0.0427 & Pred. Err. & 0.0734\\
\hline
\end{tabular}
\captionof{figure}{Accuracy of the predicted solutions to circular PEC with varying frequencies: limited testing data.}
\label{tbl:cfreqexin}
\end{center}

We next change the circular PEC to a square centered at $(0.5,0.5)$ with the side length of $0.4$. We demonstrate that the network can still predict high-fidelity solutions despite the corner singularities of the analytic solutions. We use the same training and testing data of perturbed frequencies as in the circular case. The results are displayed in Figures~\ref{fig:sq} and \ref{fig:sqc}.  
We observe from Table~\ref{tbl:sfreq} that the results are similar to those of the circular case, although the errors for each of the training gaps are more pronounced.

\begin{center}
\includegraphics[scale=0.47]{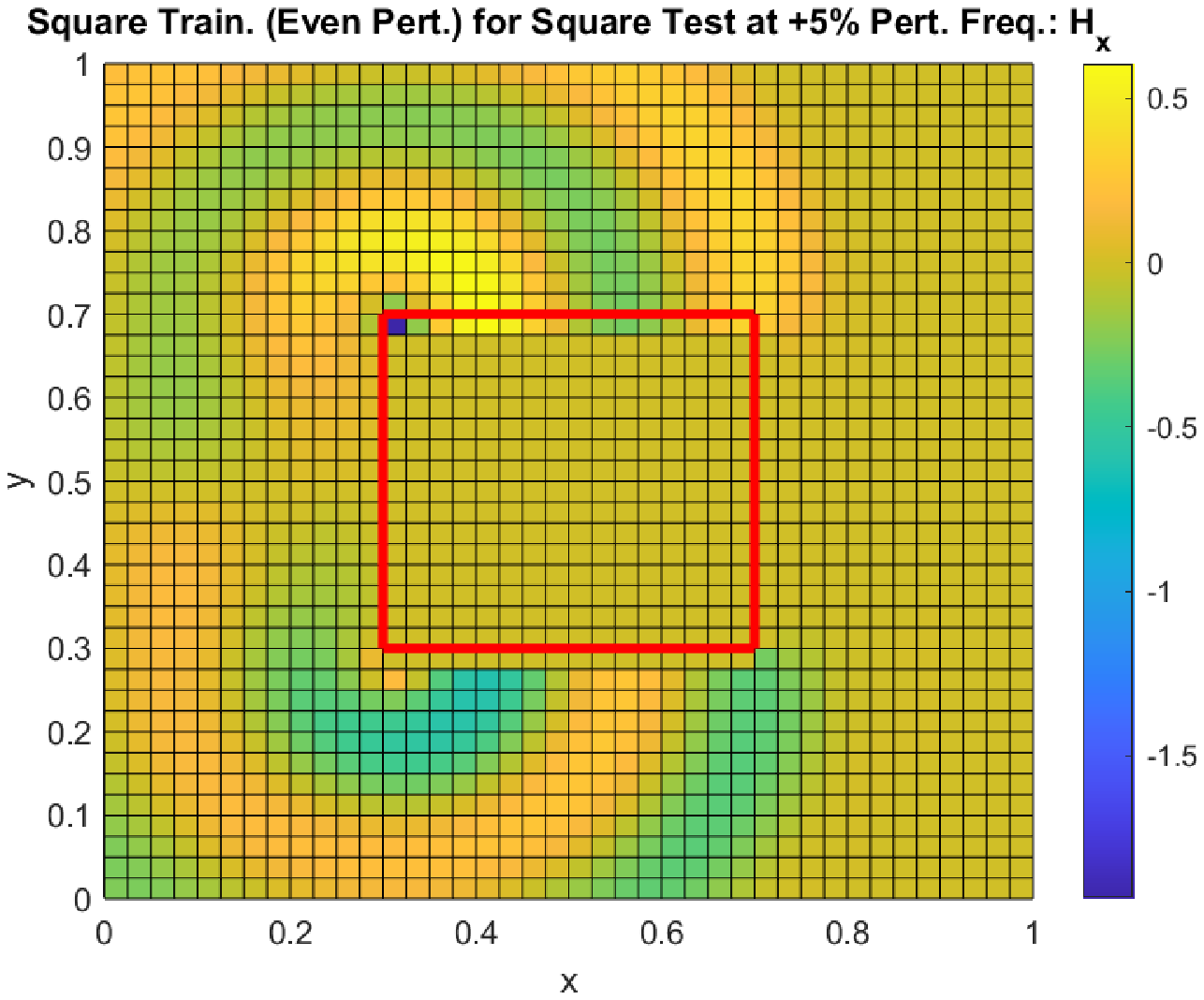}
\captionof{figure}{Predicted values of $H_x$ scattering off a square for the $+5\%$ perturbed frequency in the plane at 0.8 seconds after training the model on all the even (frequency) perturbations.} \label{fig:sq}
\end{center}

\begin{center}
\includegraphics[scale=0.47]{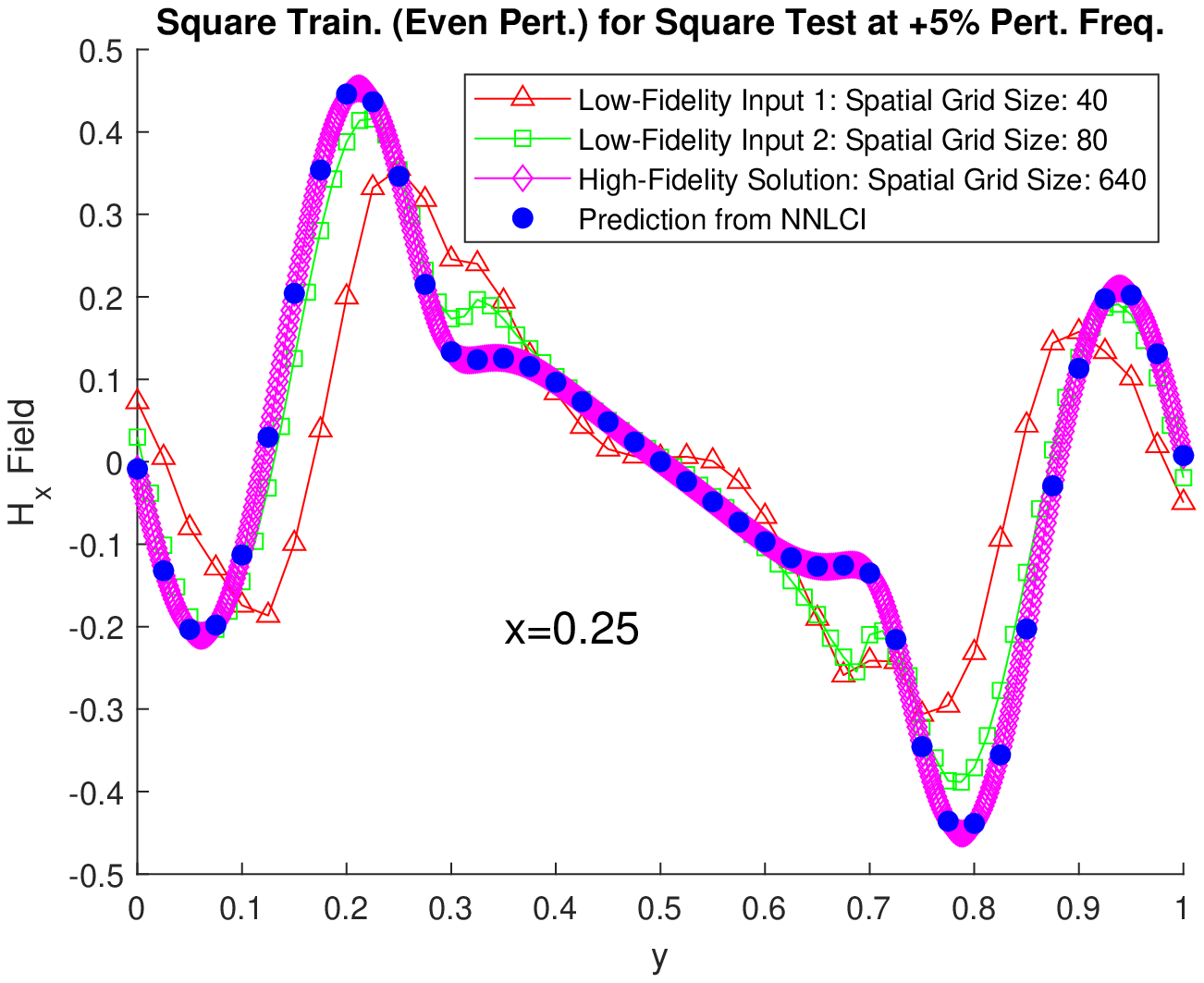}
\includegraphics[scale=0.47]{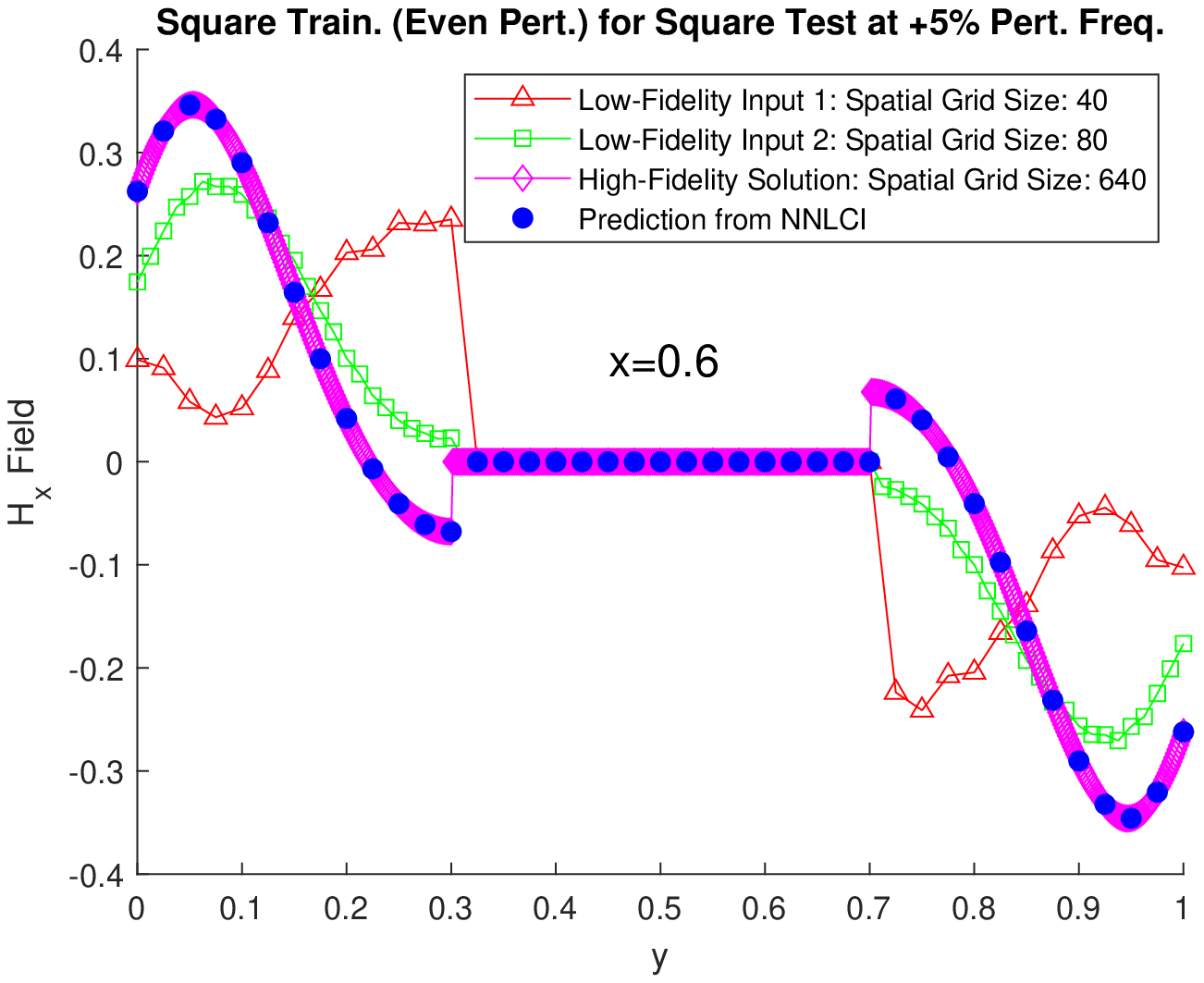}
\captionof{figure}{Cross-sections of the computed $H_x$ at different $y$ values.} \label{fig:sqc}
\end{center}

\begin{center}
\begin{tabular}{ |c|c|c|c|c|c|  }
 \hline
 \multicolumn{6}{|c|}{Square PEC: modifying frequency} \\
 \hline
 \multicolumn{2}{|c|}{Even Perturbations} & \multicolumn{2}{|c|}{$0\%,\pm 10\%$} & \multicolumn{2}{|c|}{$0\%,\pm 5\%$}\\
 \hline
Finer Input Err. & 0.1500 & Finer Input Err. & 0.1501 & Finer Input Err. & 0.1505\\
Pred. Err. & 0.0095 & Pred. Err. & 0.1006 & Pred. Err. & 0.0761\\
 \hline
\end{tabular}
\captionof{figure}{Accuracy of the predicted solutions to square PEC with varying frequencies.}\label{tbl:sfreq}
\end{center}

\end{subsubsection}

\begin{subsubsection}{Rotations}
It is well known that the Maxwell's equations \eqref{eqn:mw} are rotation invariant.
If the PEC object is also symmetric with respect to rotation as in the circular case, then the solutions are simply rotated when the planar incident waves are rotated. However, we observe poor accuracy improvement when the network is trained on one orientation of the incident waves and tested on another configuration. This can be explained by the fact that the training is entirely local, hence the network is not expected to discern the global rotation invariance of the training and testing data. On the other hand one can easily enrich the training data set with the rotated numerical solutions that do not need to be re-computed. Moreover many practical applications involve general shapes of PECs for which the rotational symmetry of the solutions is absent or unknown a priori. Hence the training of the neural networks is expected to involve incident waves entering the domain at various incident angles, regardless of whether the training itself is local or global.

\end{subsubsection}
\end{subsection}

\end{section}

\begin{subsection}{Varying PEC boundaries}
We now consider applications of our neural network when the training and testing data involve the PECs of different shapes. The complexity of the associated geometric parameter space is more complicated than that of the perturbed incident waves. The local training of the network is expected to play a crucial role of rendering model order reduction to variations in local boundary geometries of the PECs. Unless stated otherwise the incident waves are the unperturbed ones prescribed in \eqref{eqn:ww}.

\begin{subsubsection}{Modification of the circular PEC}
We first consider varying the radius of the circular PEC. We generate the data corresponding to the perturbations of the unperturbed radius of $0.2$ with perturbations of $0,\pm 2\%, \pm 3\%, \pm 4\%,$ $\pm 5\%, \pm 6\%, \pm 7\%, \pm 8\%, \pm 10\%$, just as for the previously considered cases. Similar to the case of varying amplitudes, we separate the generated data into either training or testing data. We observe significant accuracy improvement even when the number of training data is quite small.
\begin{center}
\begin{tabular}{ |c|c|c|c|c|c|  }
 \hline
 \multicolumn{6}{|c|}{Circular PEC: modifying radius} \\
 \hline
 \multicolumn{2}{|c|}{Even Perturbations} & \multicolumn{2}{|c|}{$0\%,\pm 10\%$} & \multicolumn{2}{|c|}{$0\%,\pm 5\%$}\\
 \hline
Finer Input Err. & 0.1625 & Finer Input Err. & 0.1624 & Finer Input Err. & 0.1623\\
Pred. Err. & 0.0086 & Pred. Err. & 0.0176 & Pred. Err. & 0.0163\\
 \hline
\end{tabular}
\captionof{table}{Accuracy of predicted solutions subject to changing radius}
\label{tbl:cr}
\end{center}

We next truncate the circular PEC to circular sectors which subtends the angles of various magnitudes (see Figures ~\ref{fig:sectors1} and \ref{fig:sectors2}). 
At first, we generated the solutions corresponding to the angles from 30 to 360 degrees in 30 degree increments. We found however that our neural network, regardless of our testing and training sets, was unable to produce accurate predictions due to the large training gap. On the other hand, our network manages to improve the accuracy of the low fidelity solutions when the increments in angles are $5$ degrees. More specifically we train when the increments in angles are $10$ degrees, that is $45^{\circ},55^{\circ},\ldots,125^{\circ},135^{\circ}$, and then test on the even perturbations $50^{\circ},60^{\circ},\ldots,120^{\circ},130^{\circ}$, see Figures~\ref{fig:90degree} and \ref{fig:90degree_cross}.  We then consider enlarging the training gap to $15^{\circ}$ increments since in some practical applications it may not be ideal to enforce smaller training gaps. Namely we train on $45^{\circ},60^{\circ},\ldots,120^{\circ},135^{\circ}$ and test on the rest of the perturbations to obtain the similar accuracy improvement as in the case of $5^{\circ}$ increments, see Table \ref{tbl:sec}. 

\begin{center}
\includegraphics[scale=0.25]{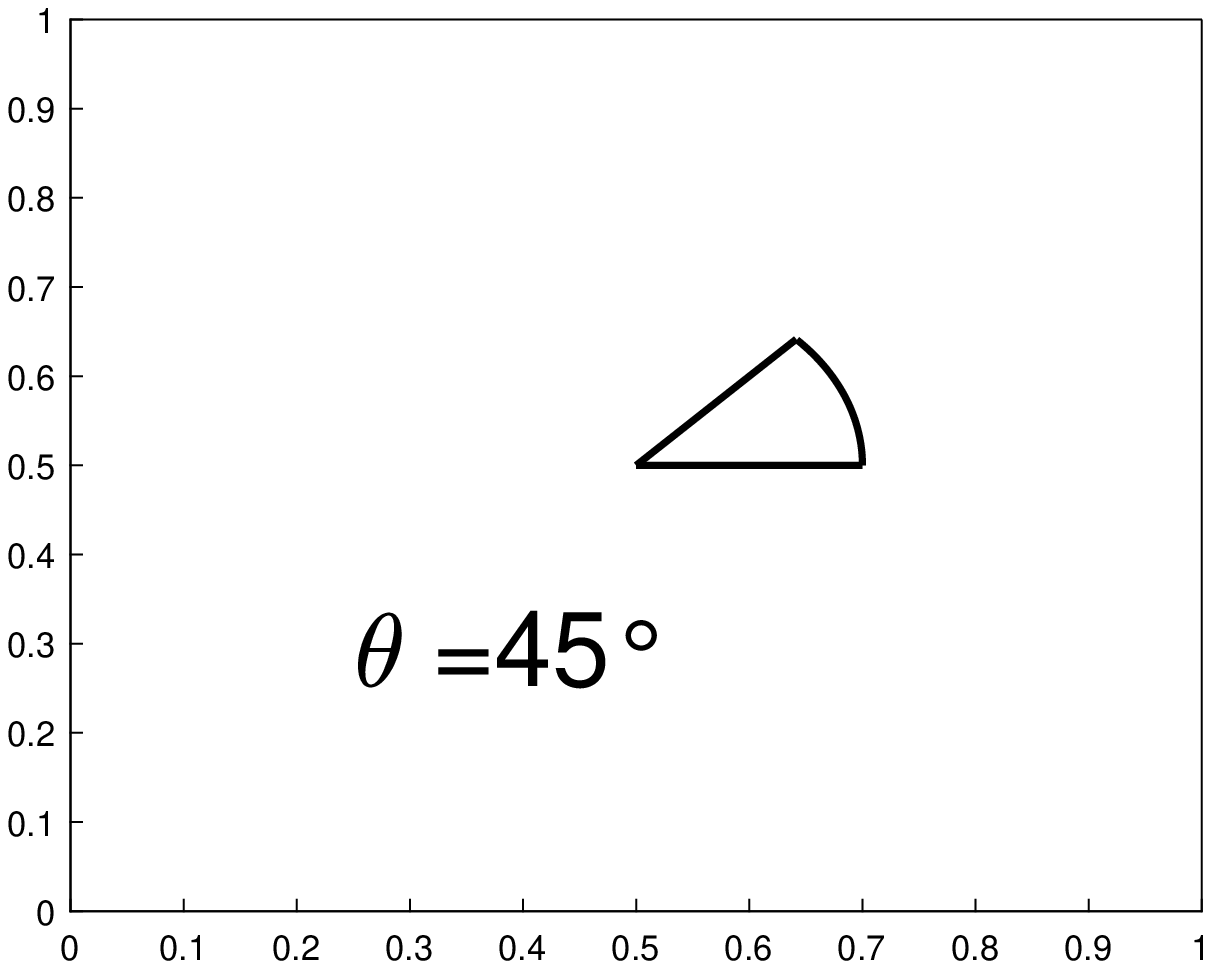}
\includegraphics[scale=0.25]{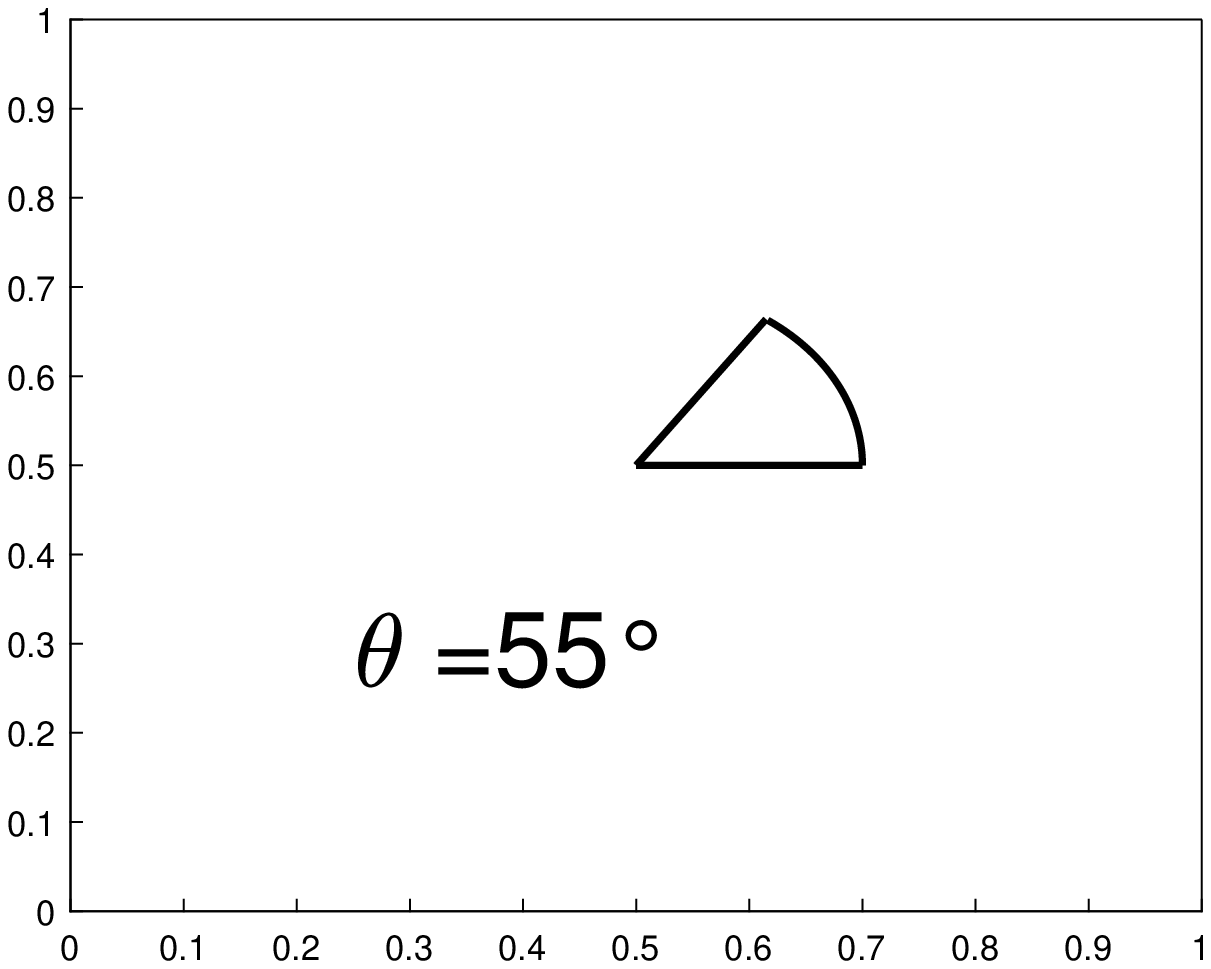}
\includegraphics[scale=0.25]{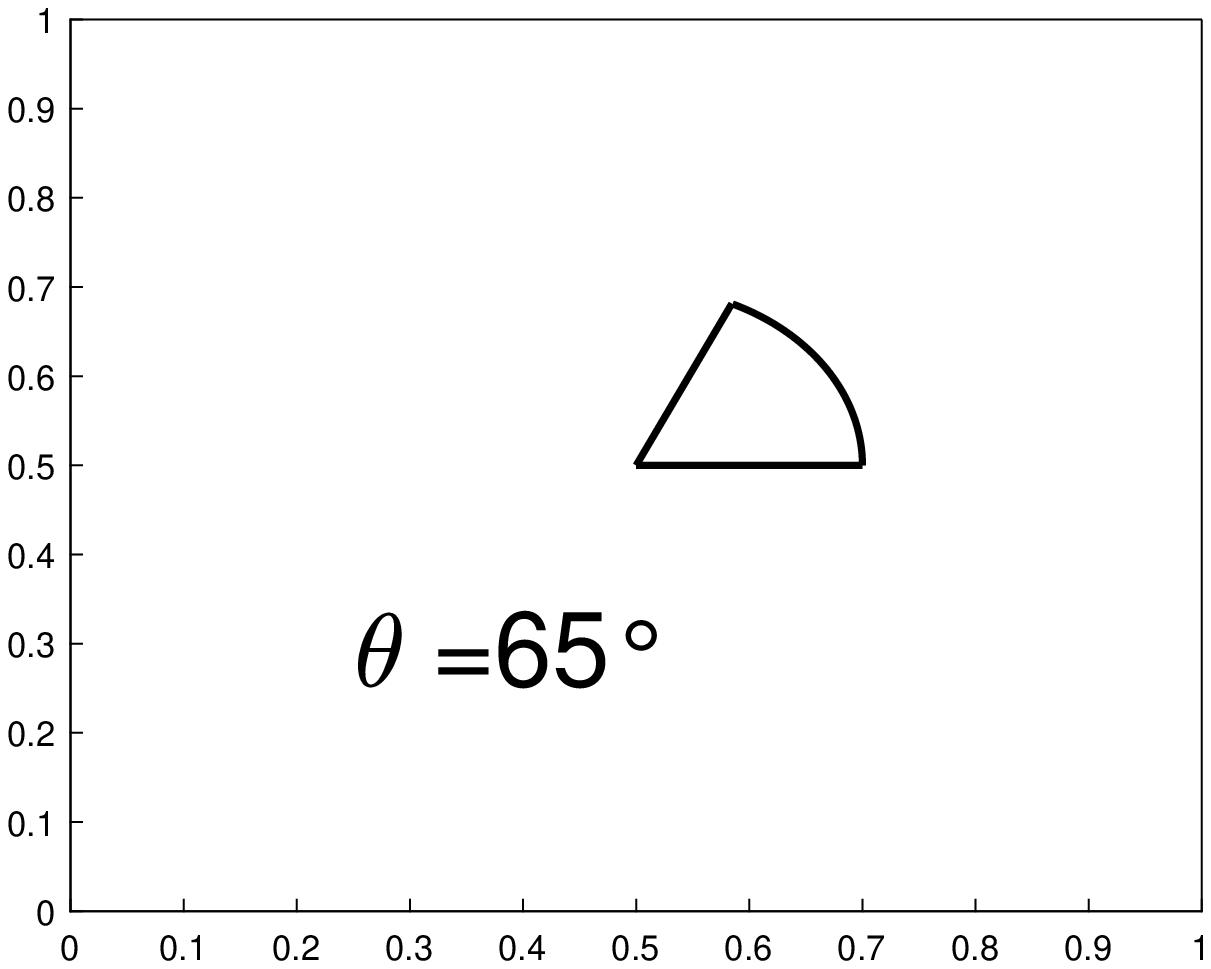}
\includegraphics[scale=0.25]{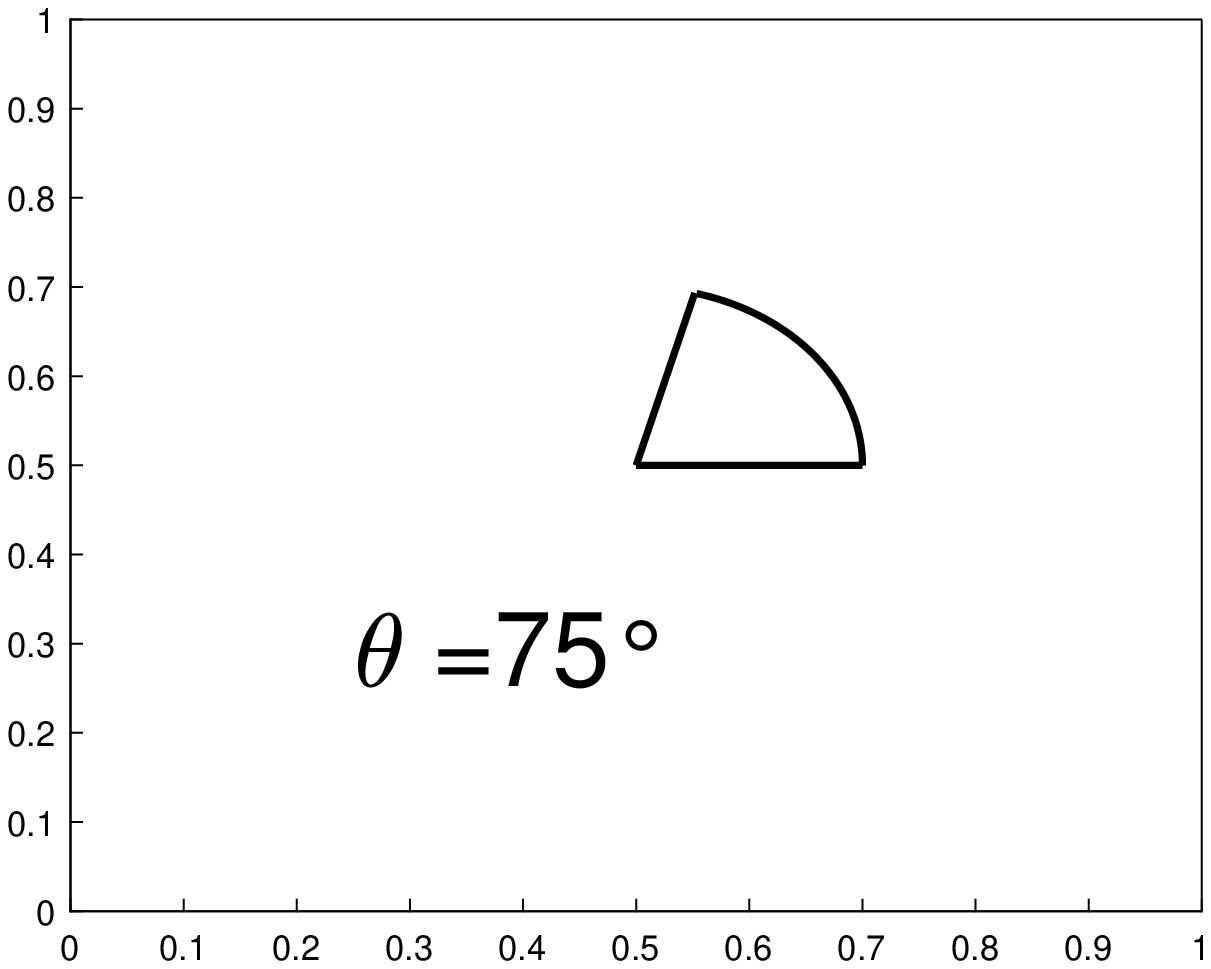}
\includegraphics[scale=0.25]{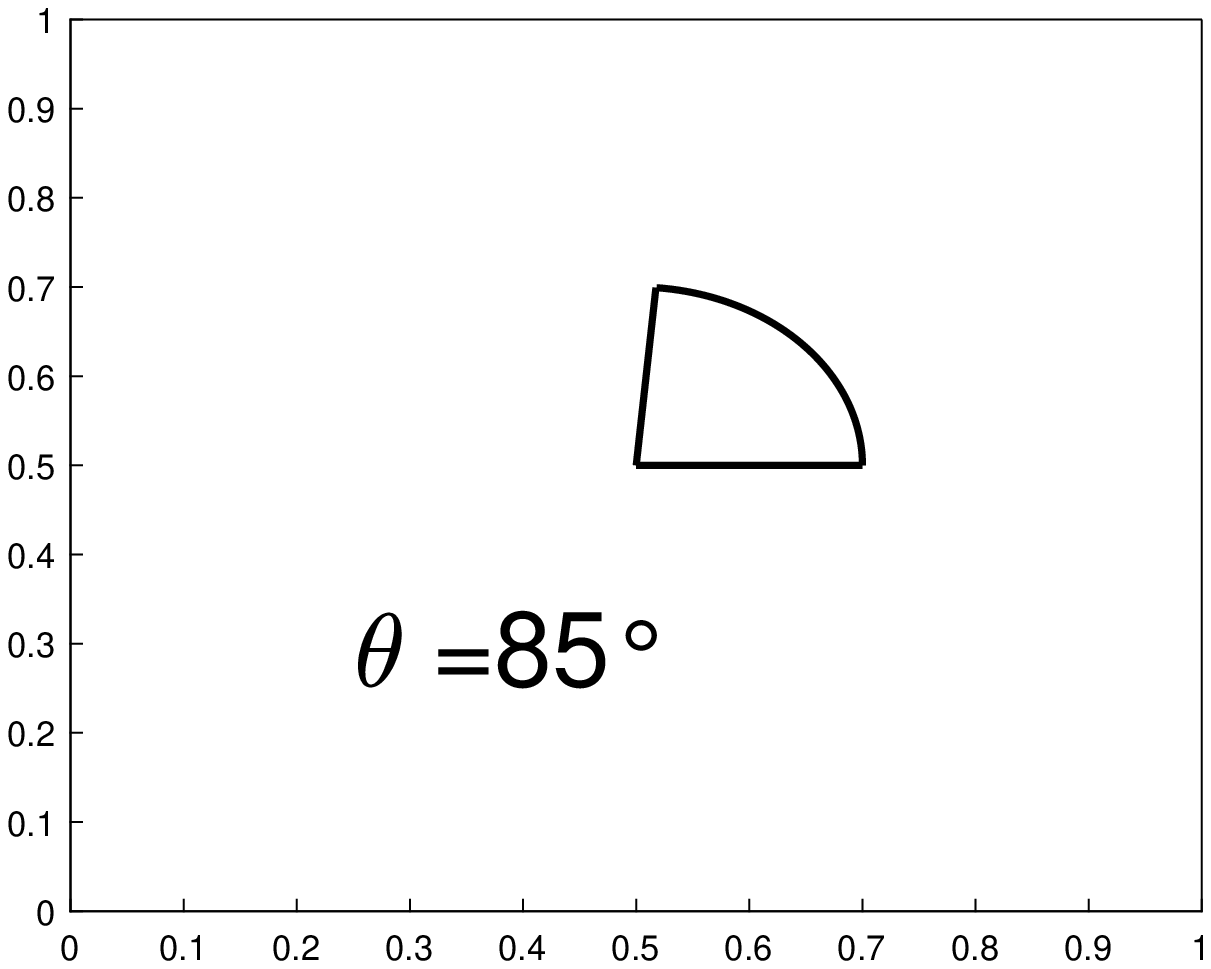}
\includegraphics[scale=0.25]{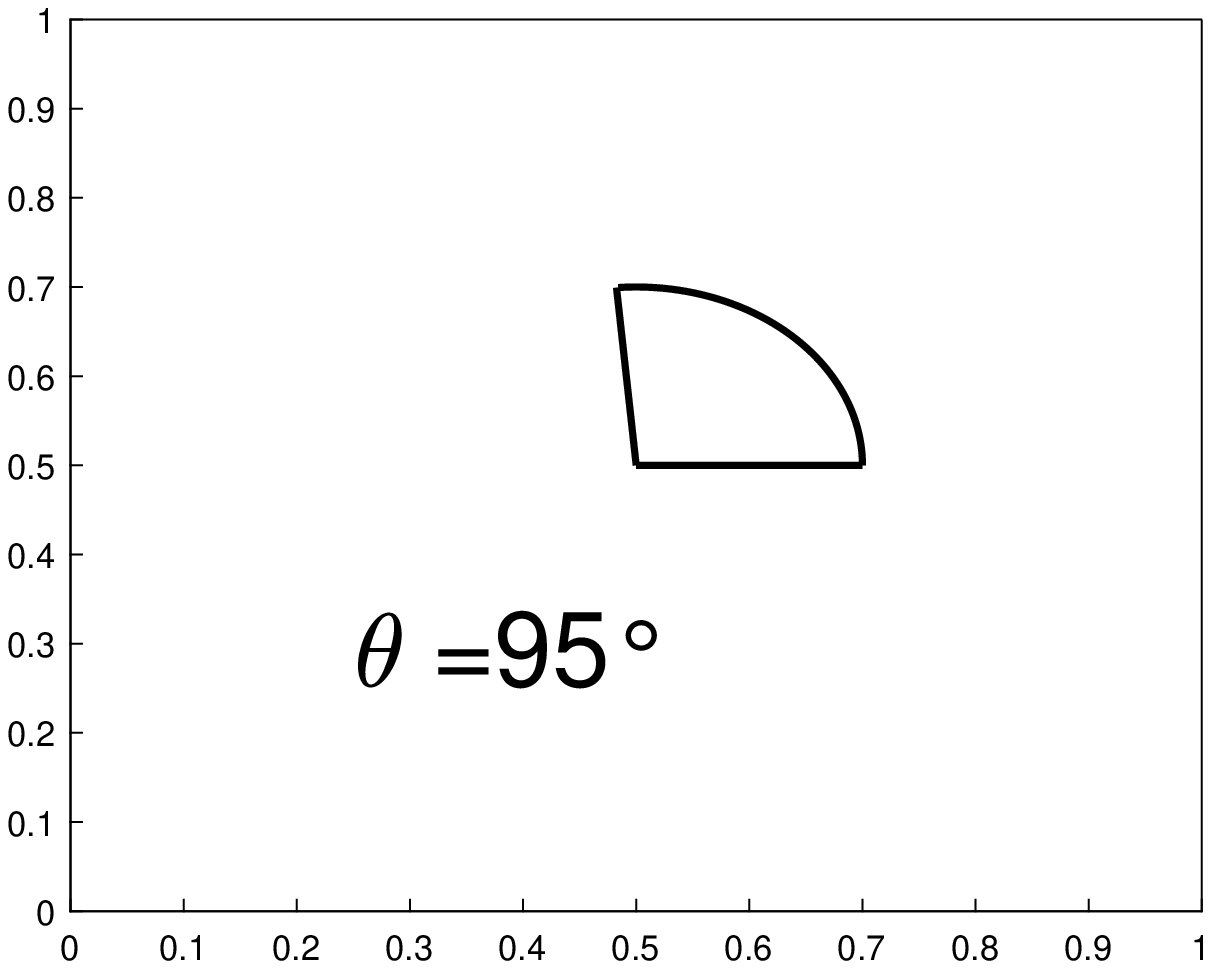}
\includegraphics[scale=0.25]{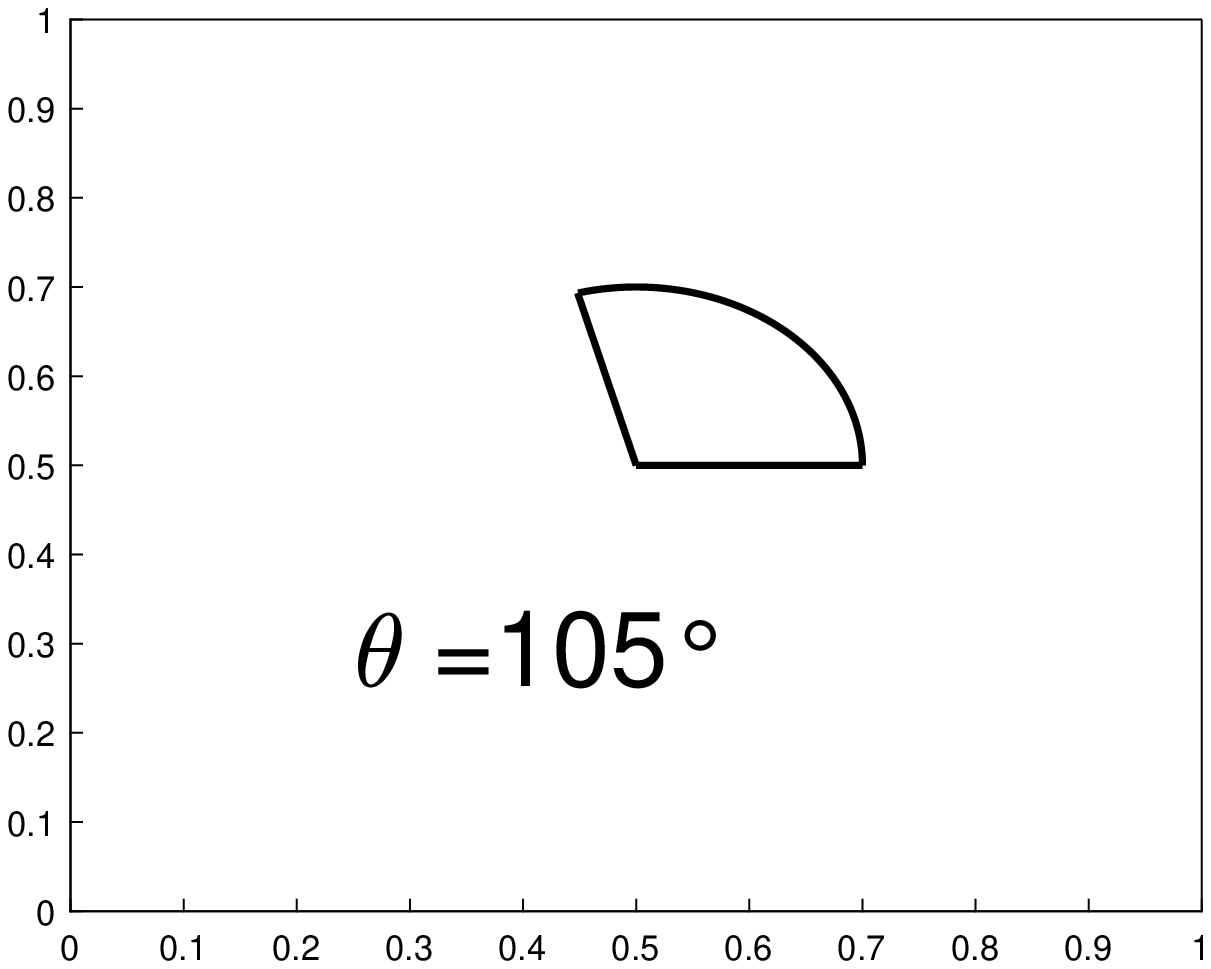}
\includegraphics[scale=0.25]{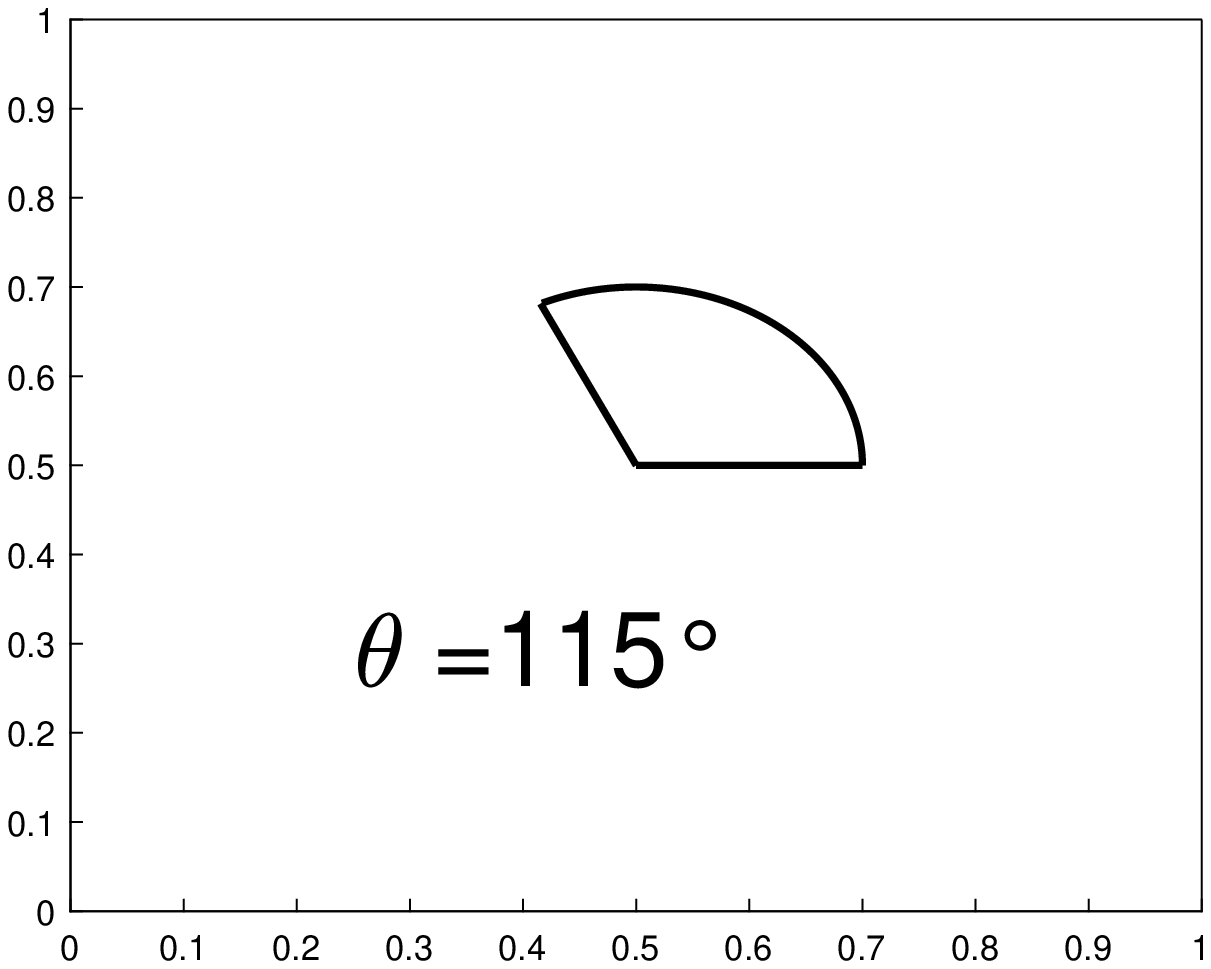}
\includegraphics[scale=0.25]{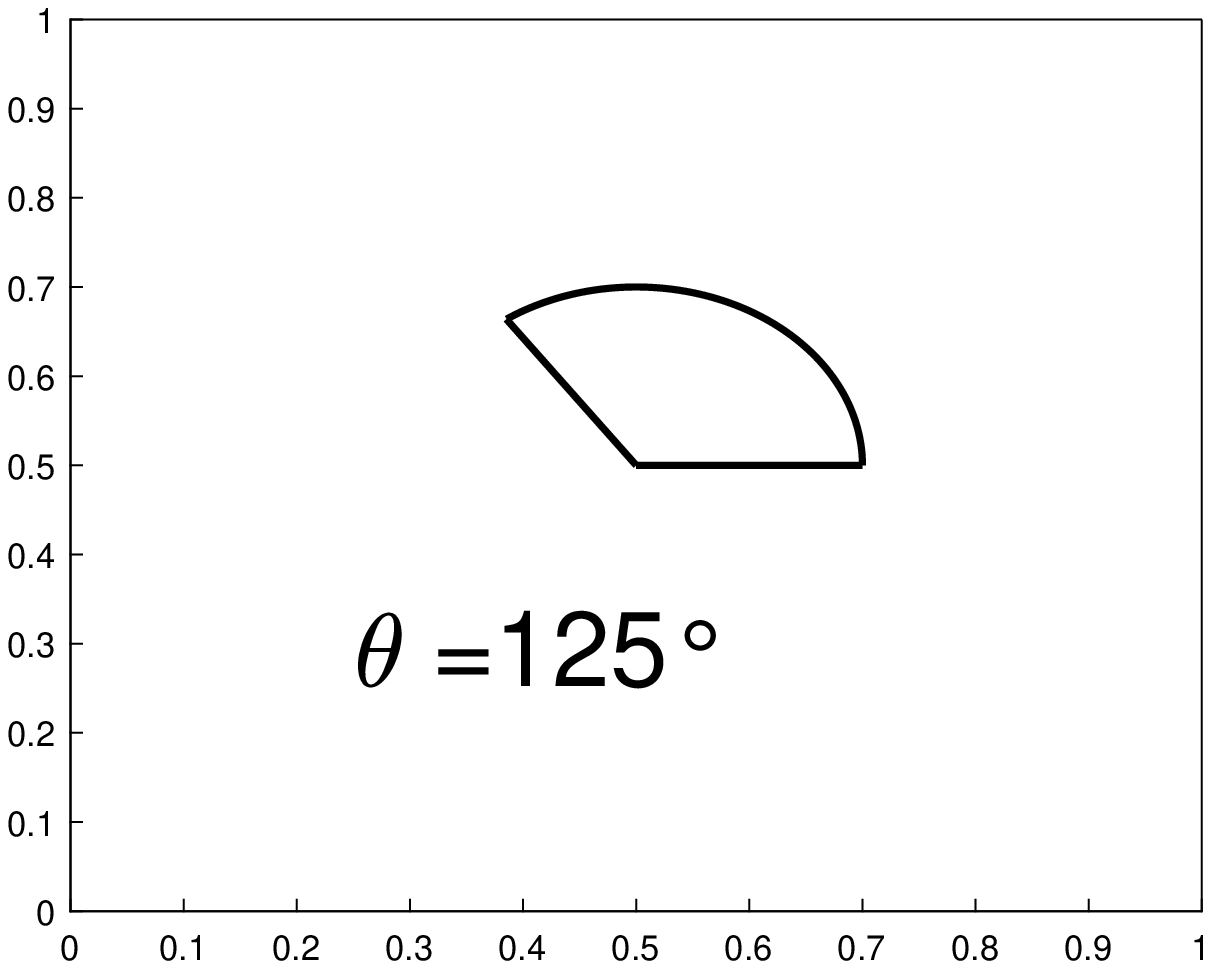}
\includegraphics[scale=0.25]{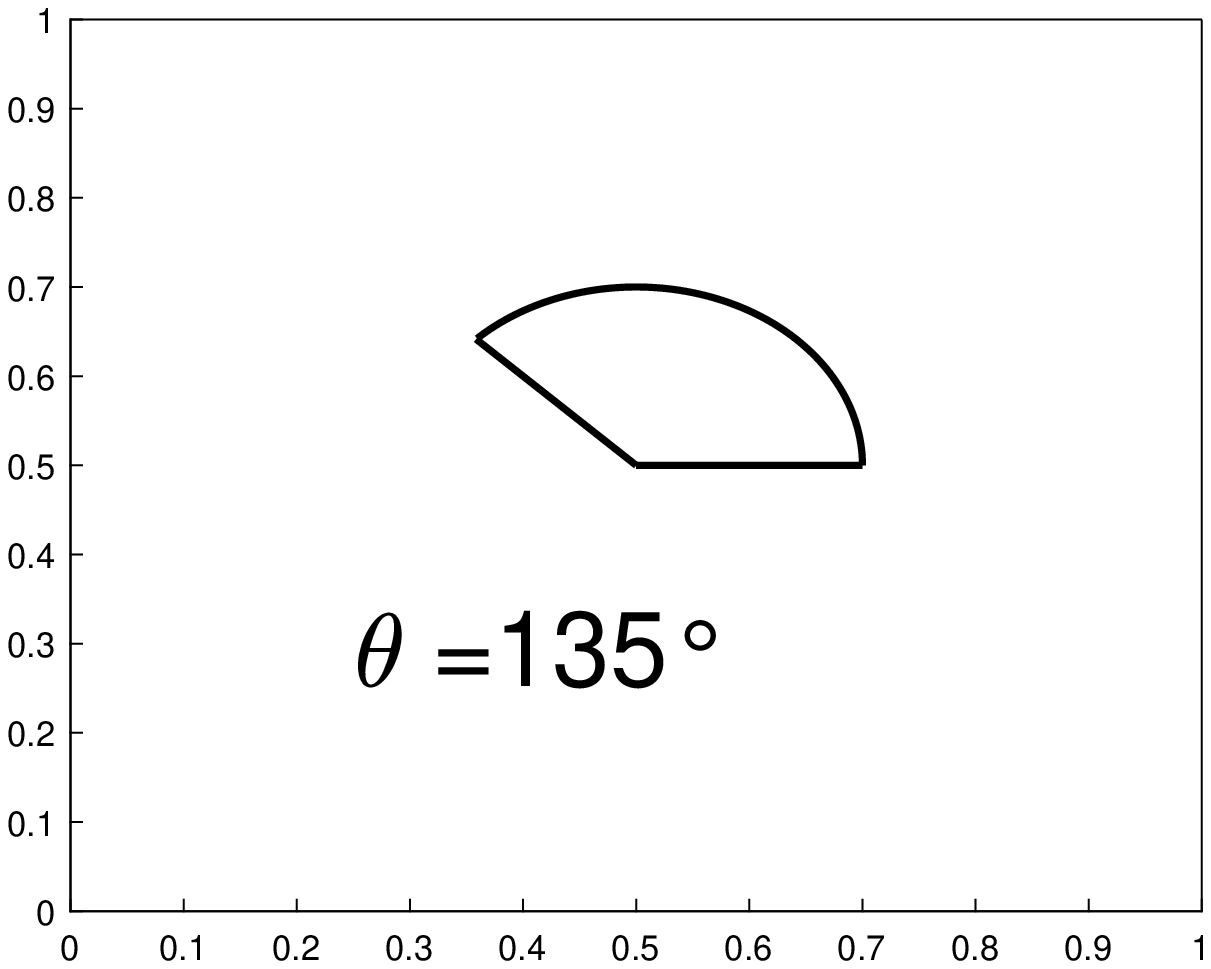}
\captionof{figure}{Shapes of Training Examples}\label{fig:sectors1}
\includegraphics[scale=0.25]{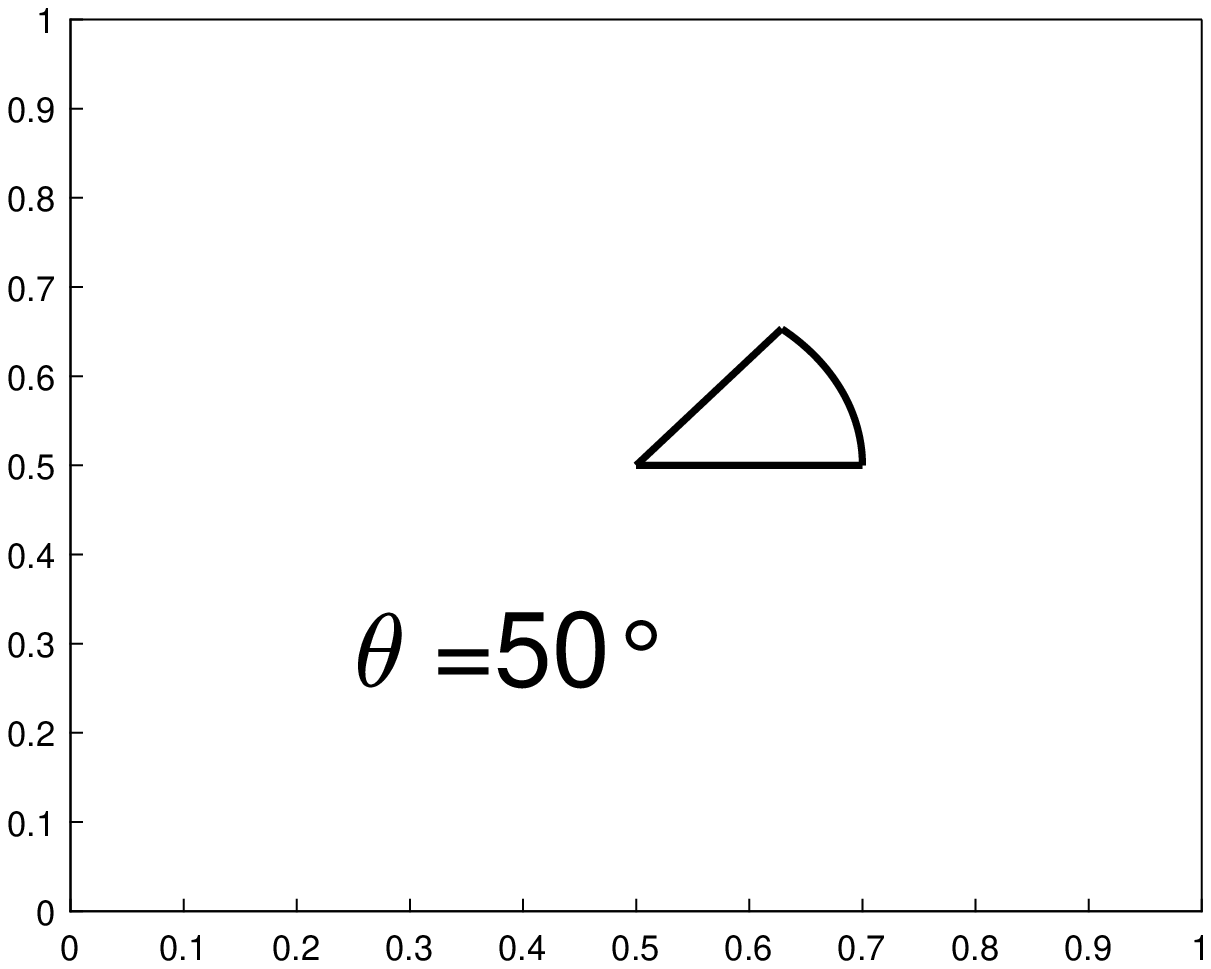}
\includegraphics[scale=0.25]{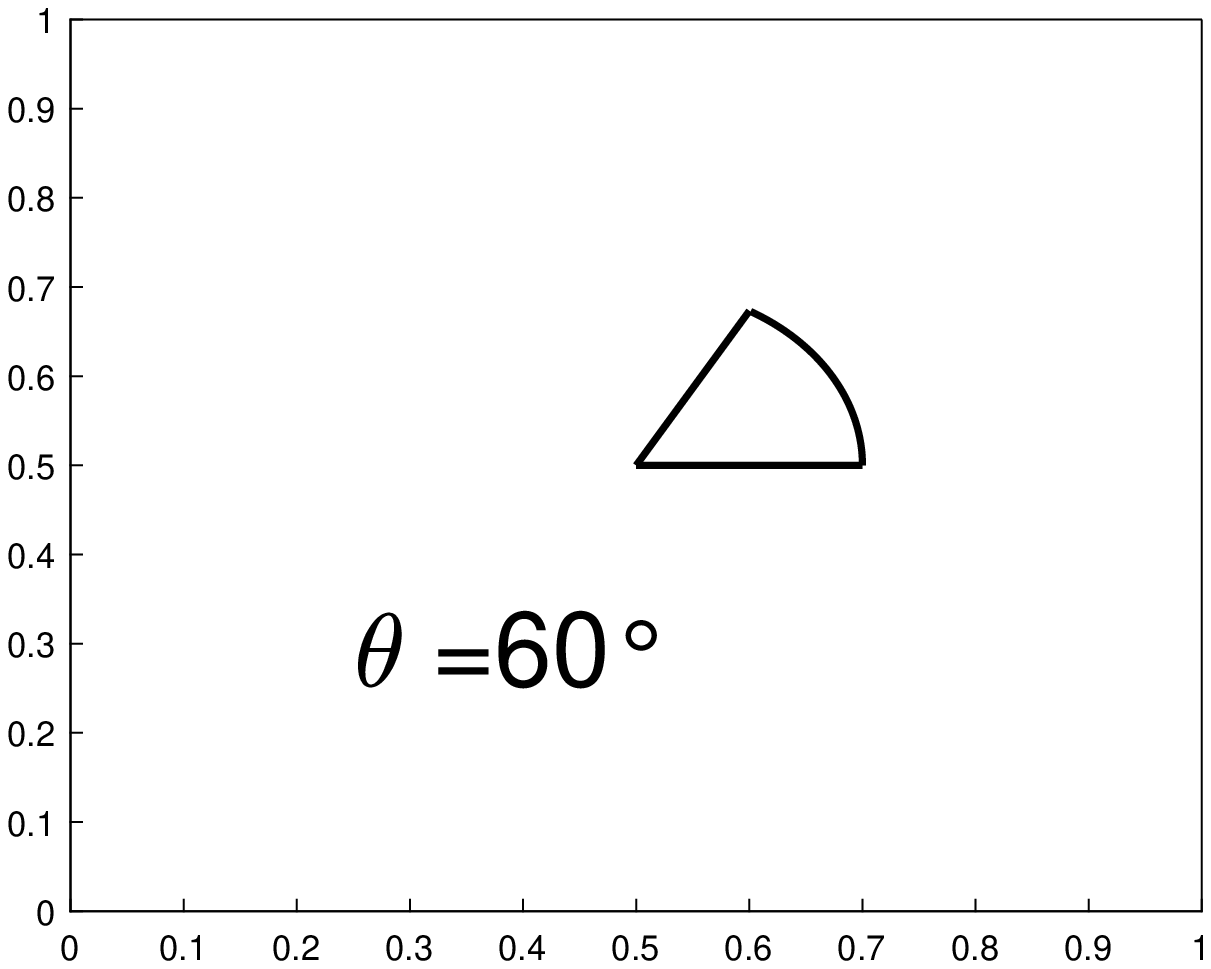}
\includegraphics[scale=0.25]{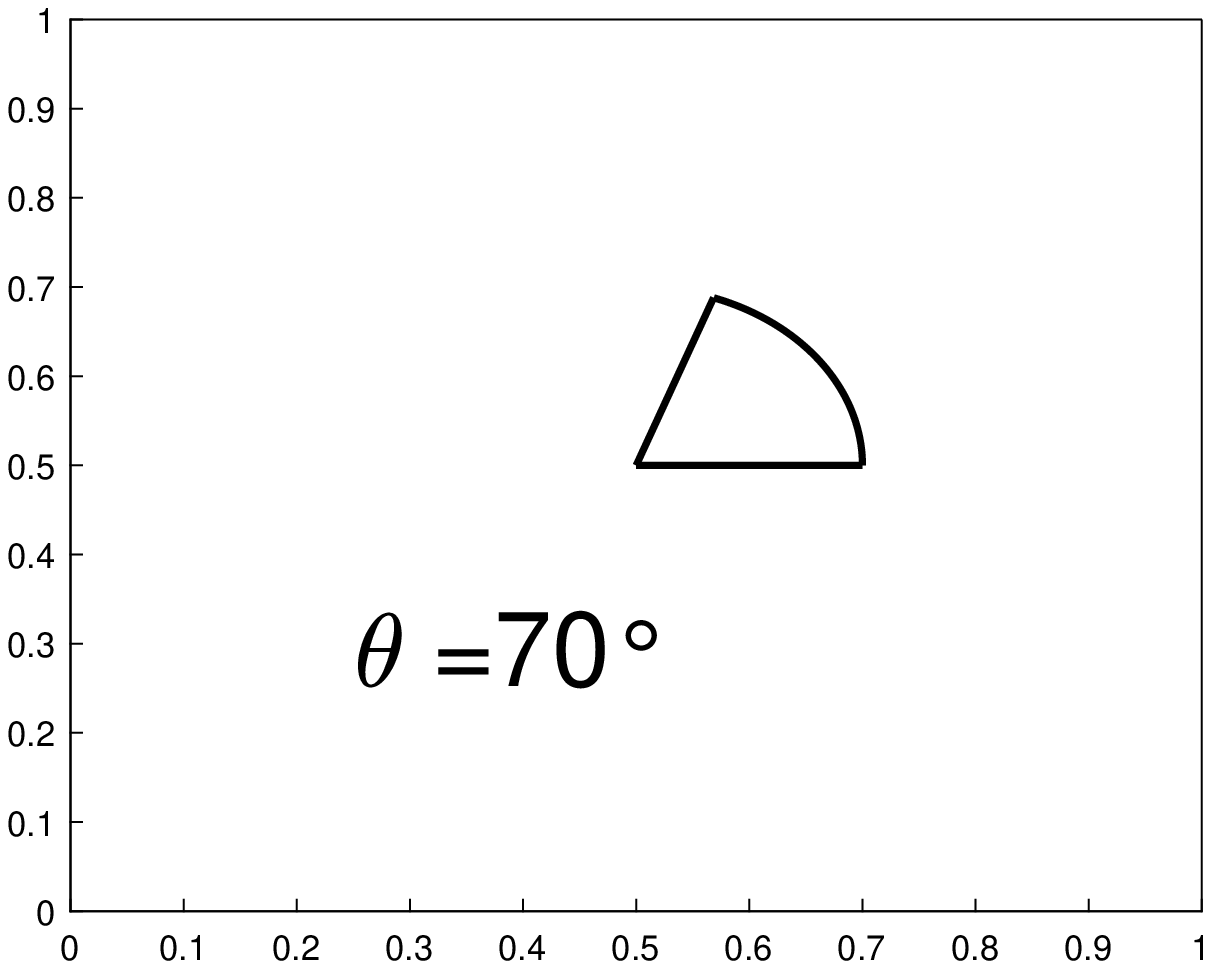}
\includegraphics[scale=0.25]{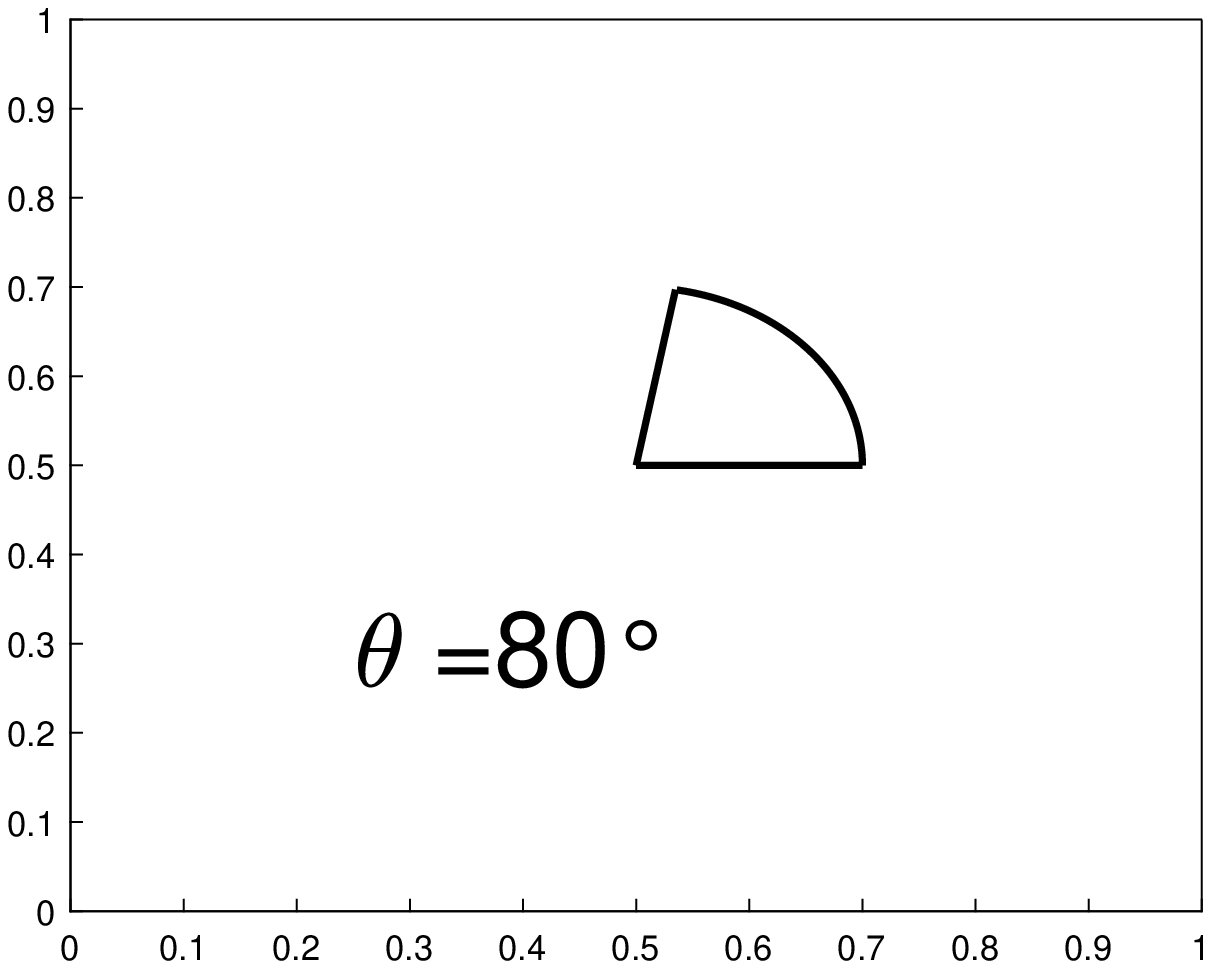}
\includegraphics[scale=0.25]{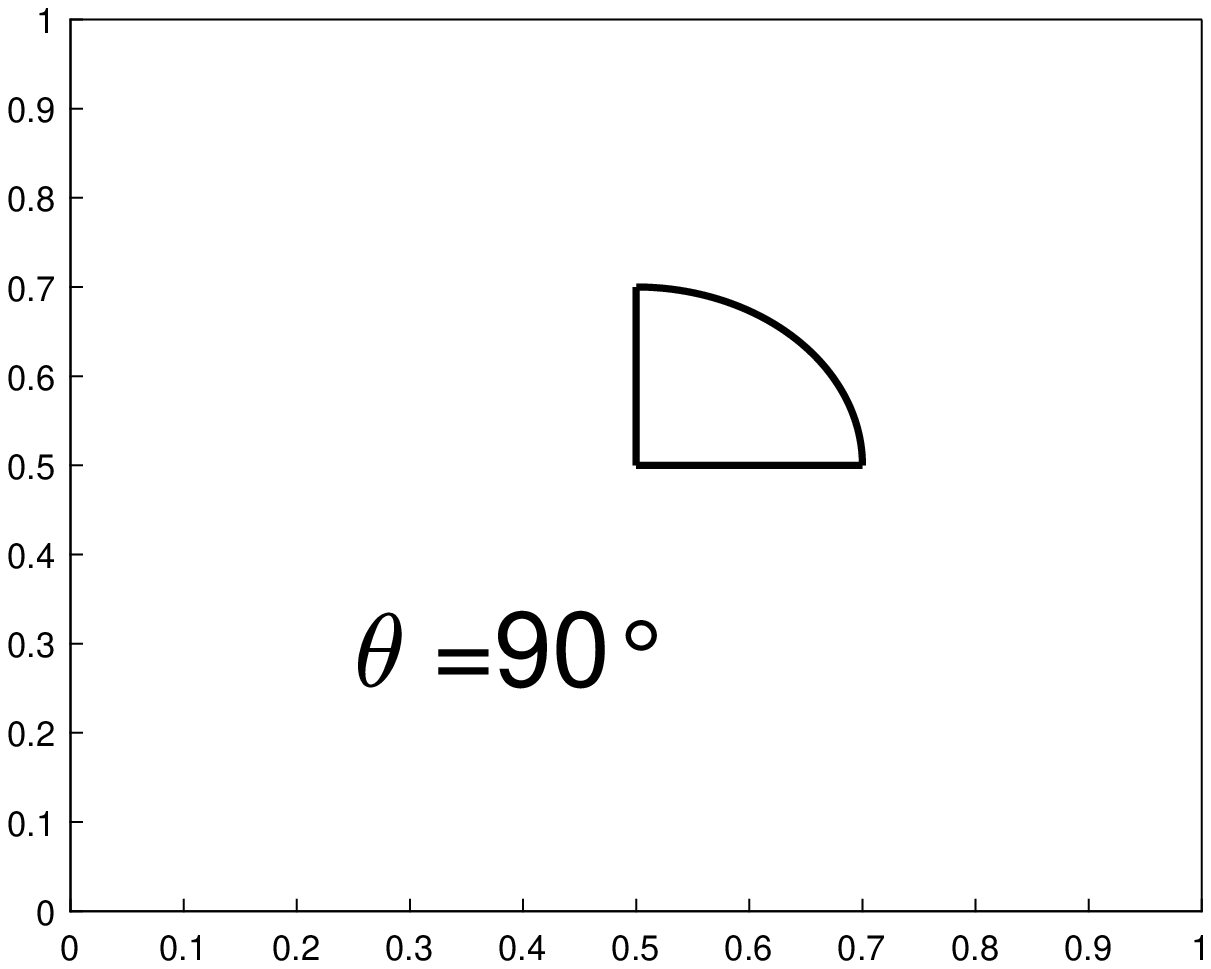}
\includegraphics[scale=0.25]{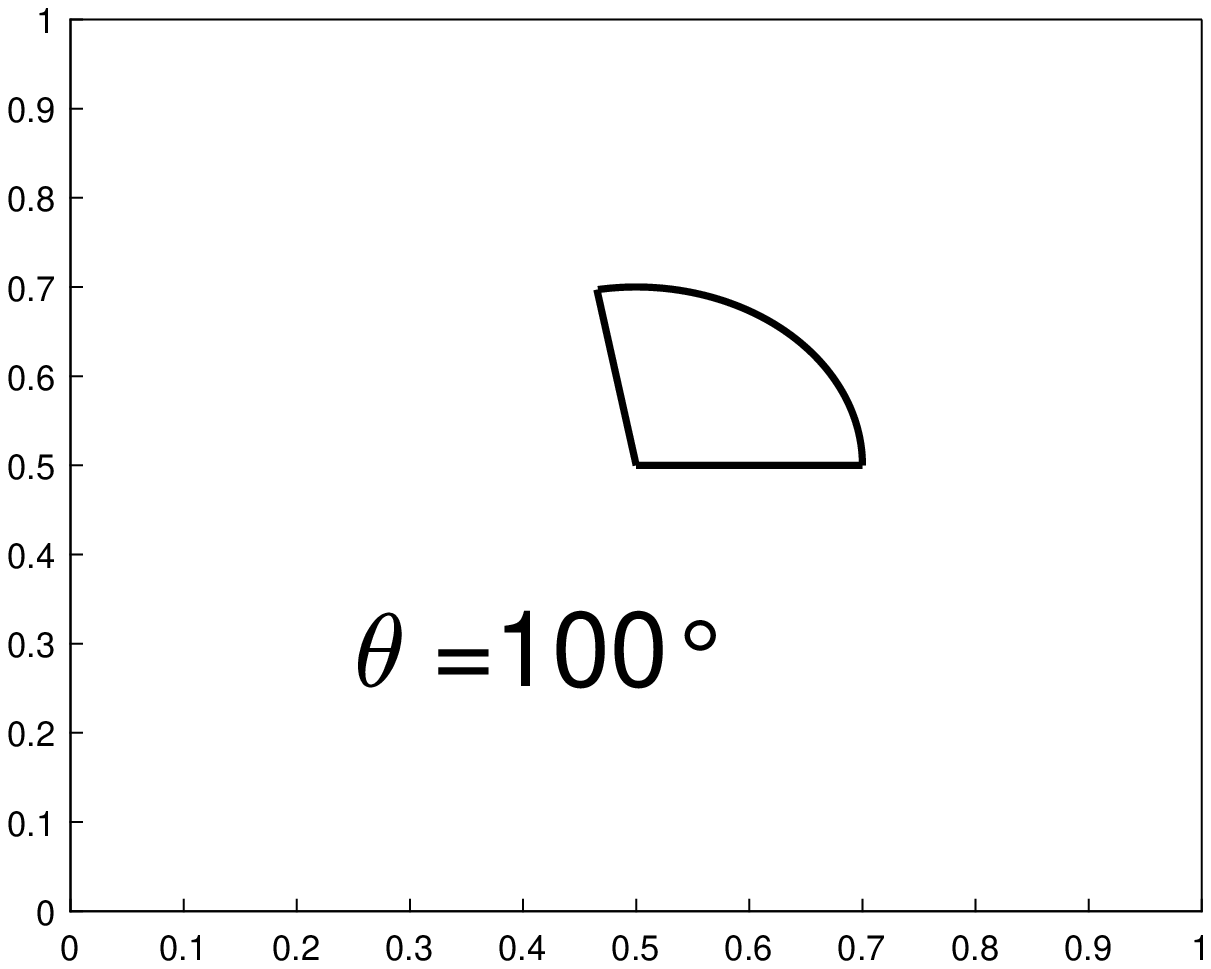}
\includegraphics[scale=0.25]{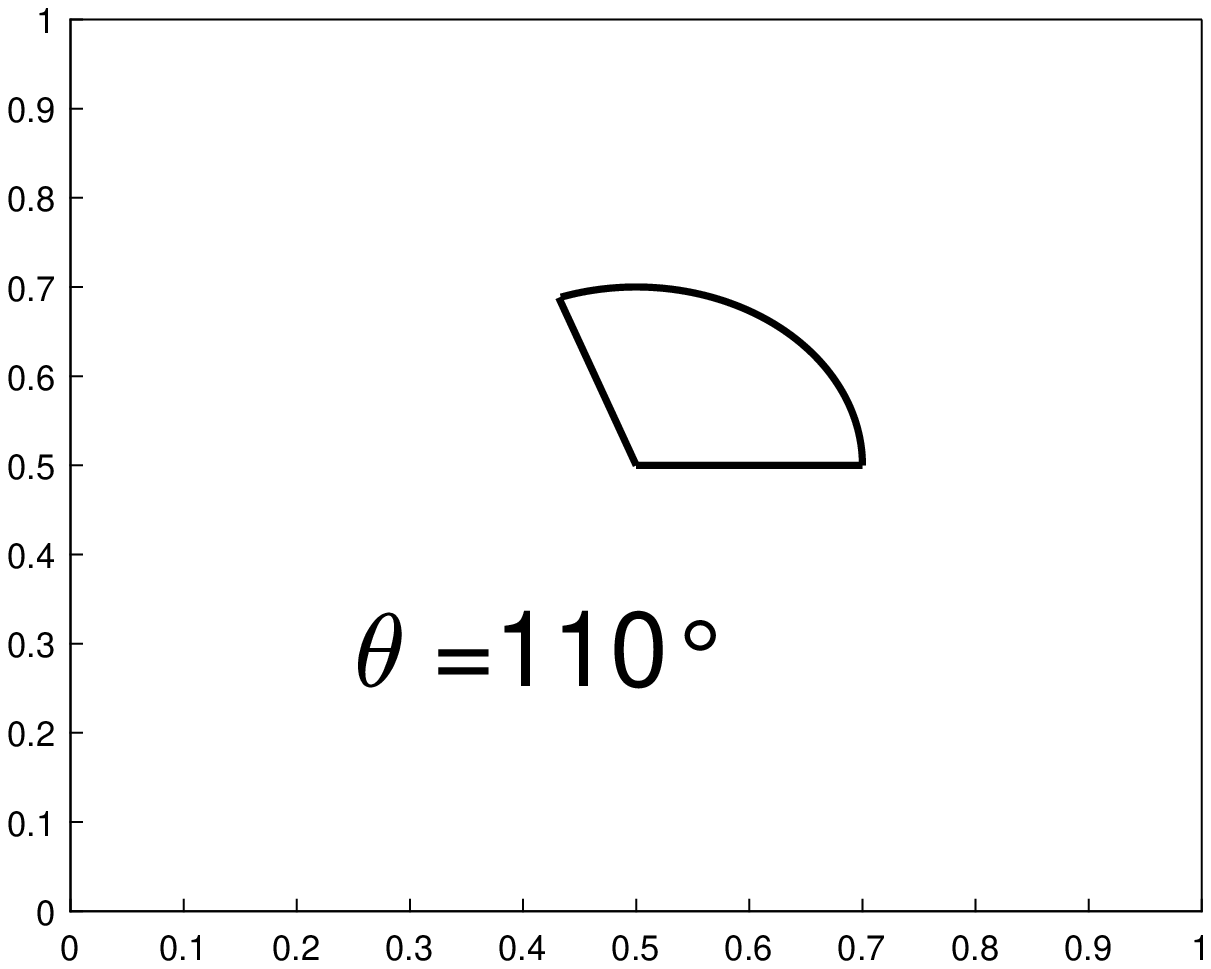}
\includegraphics[scale=0.25]{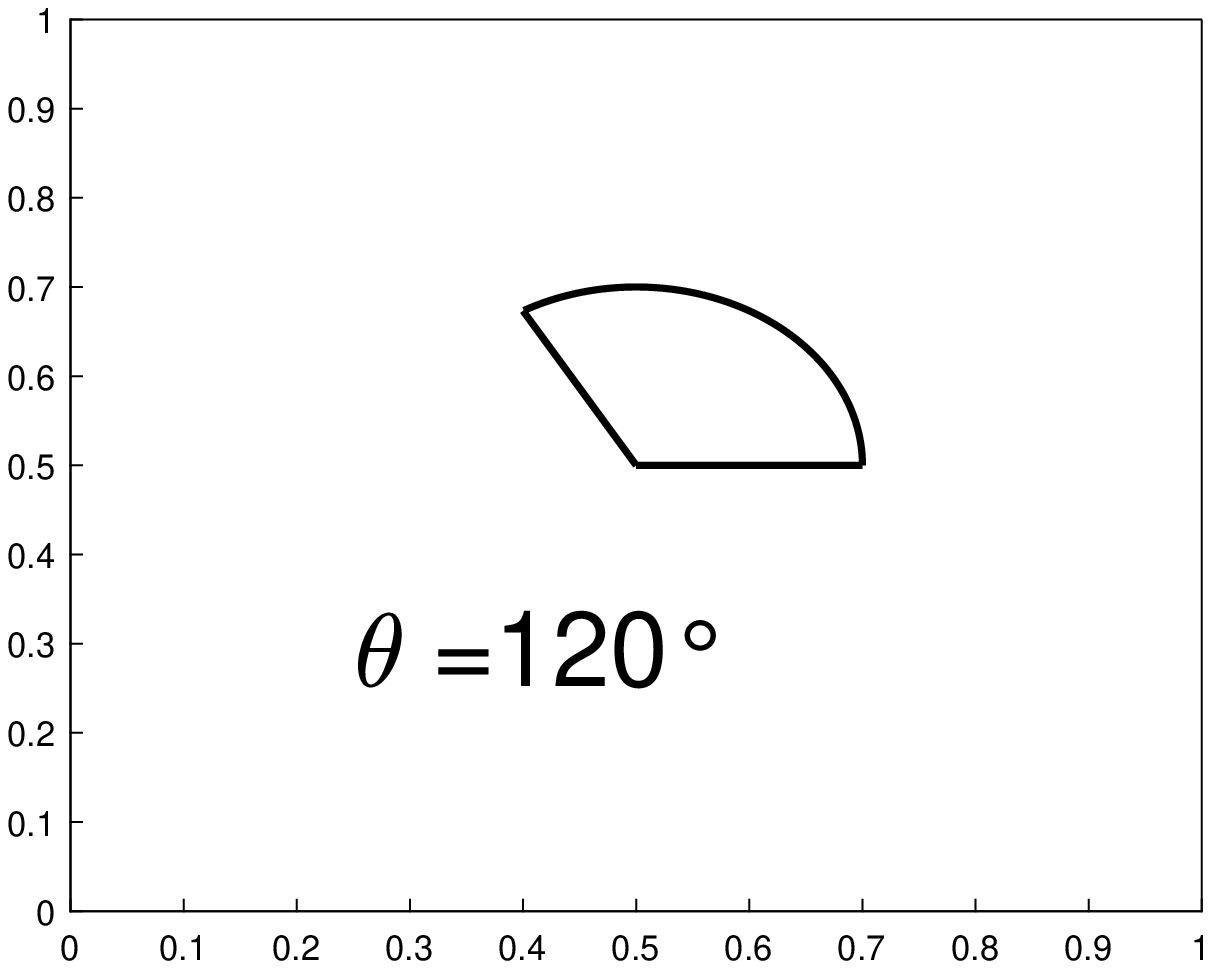}
\includegraphics[scale=0.25]{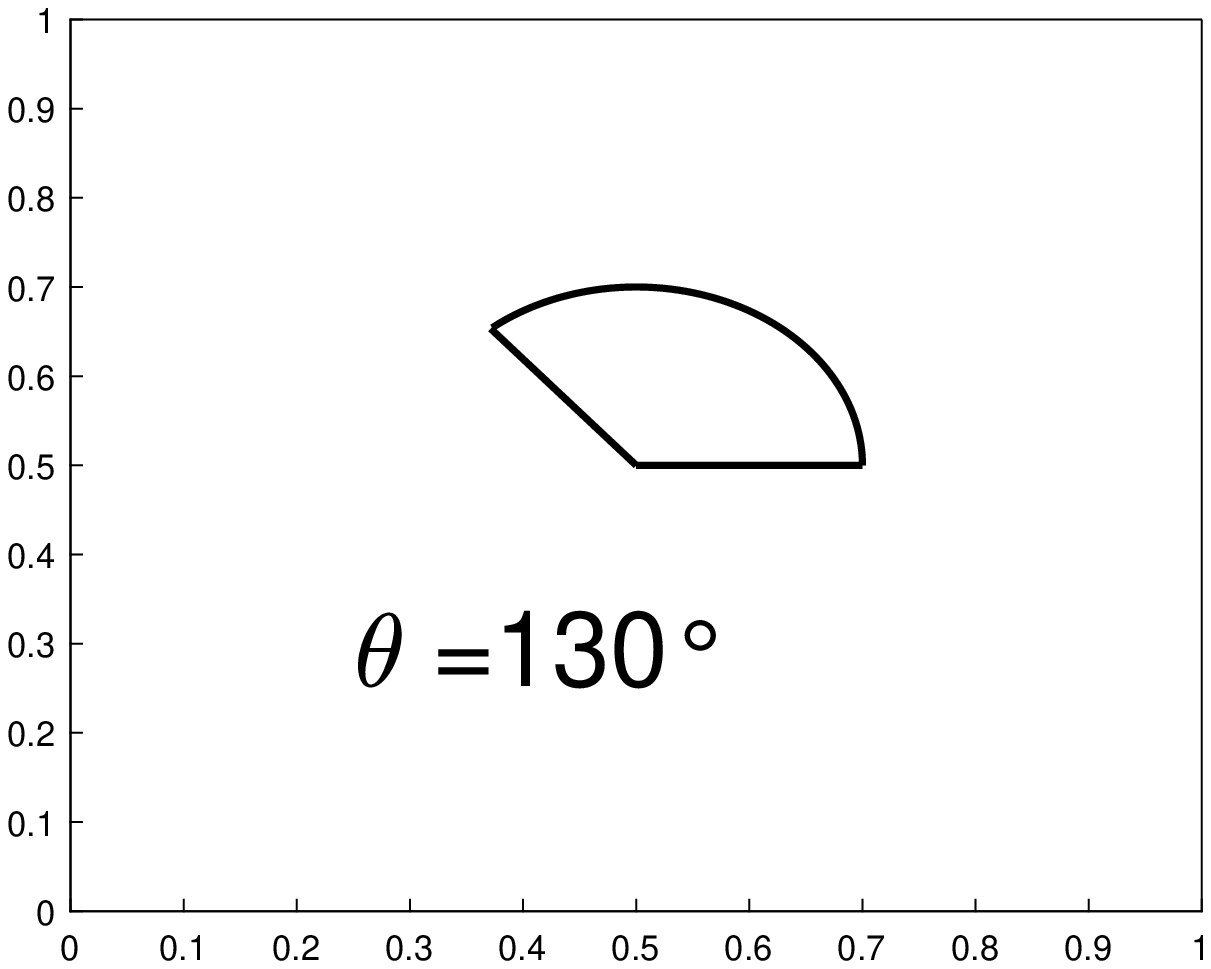}
\captionof{figure}{Shapes of Testing Examples}
\label{fig:sectors2}
\end{center}

\begin{center}
\includegraphics[scale=0.5]{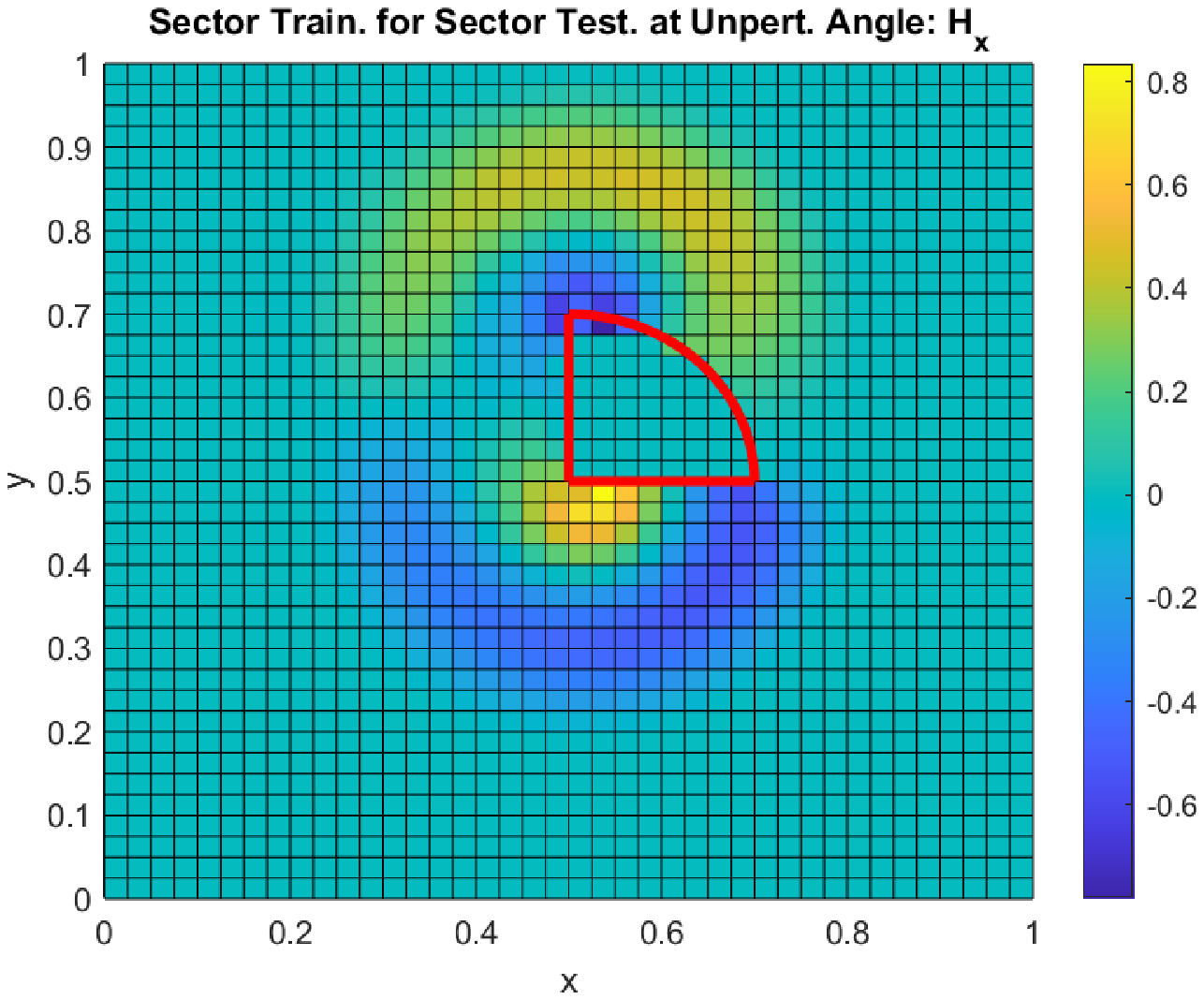}\\
\captionof{figure}{Predicted values of $H_x$ scattering off a $90^{\circ}$ sector in the plane at 0.8 seconds after training the model on  10$^\circ$ angle increments.}\label{fig:90degree}
\includegraphics[scale=0.5]{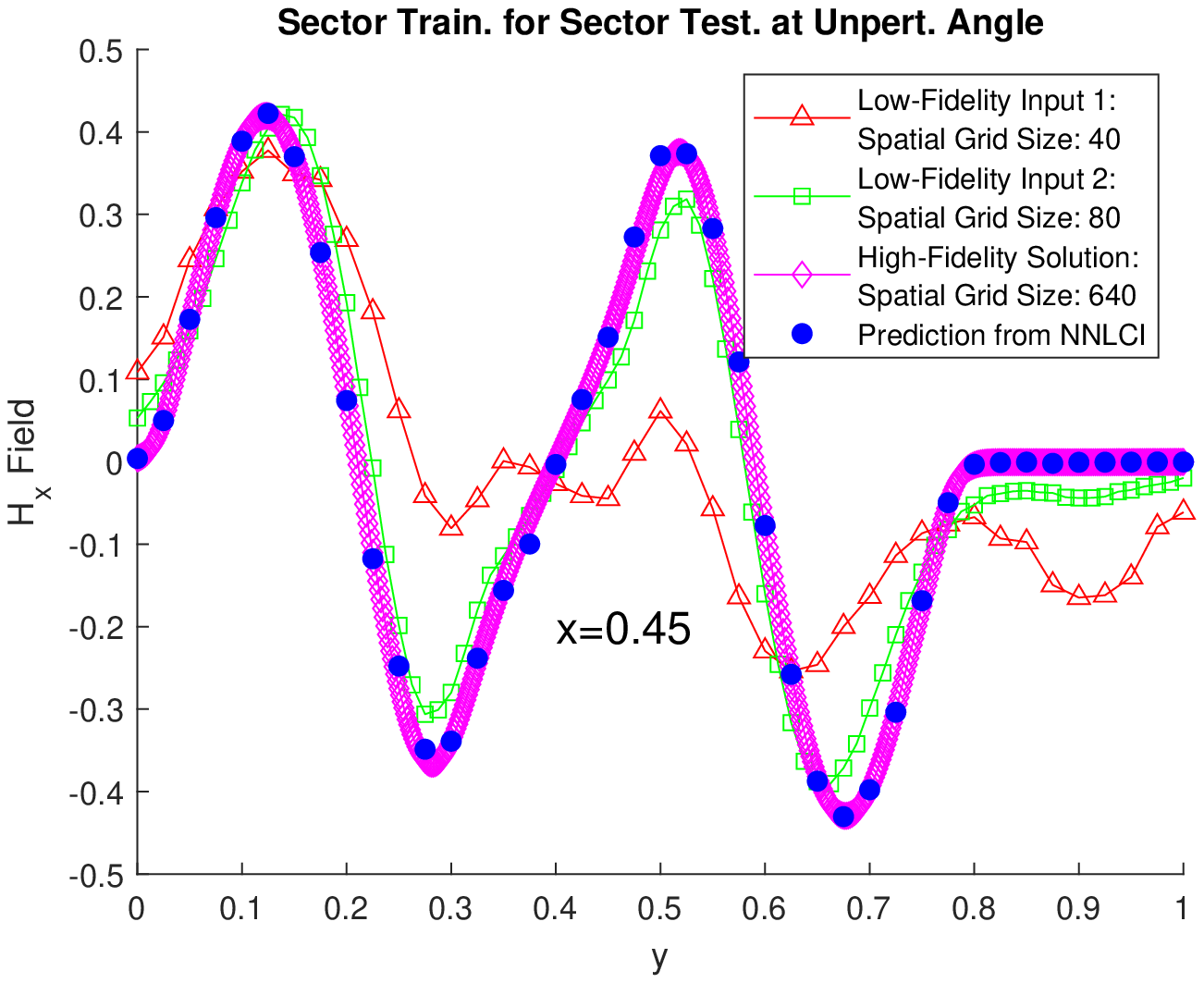}
\includegraphics[scale=0.5]{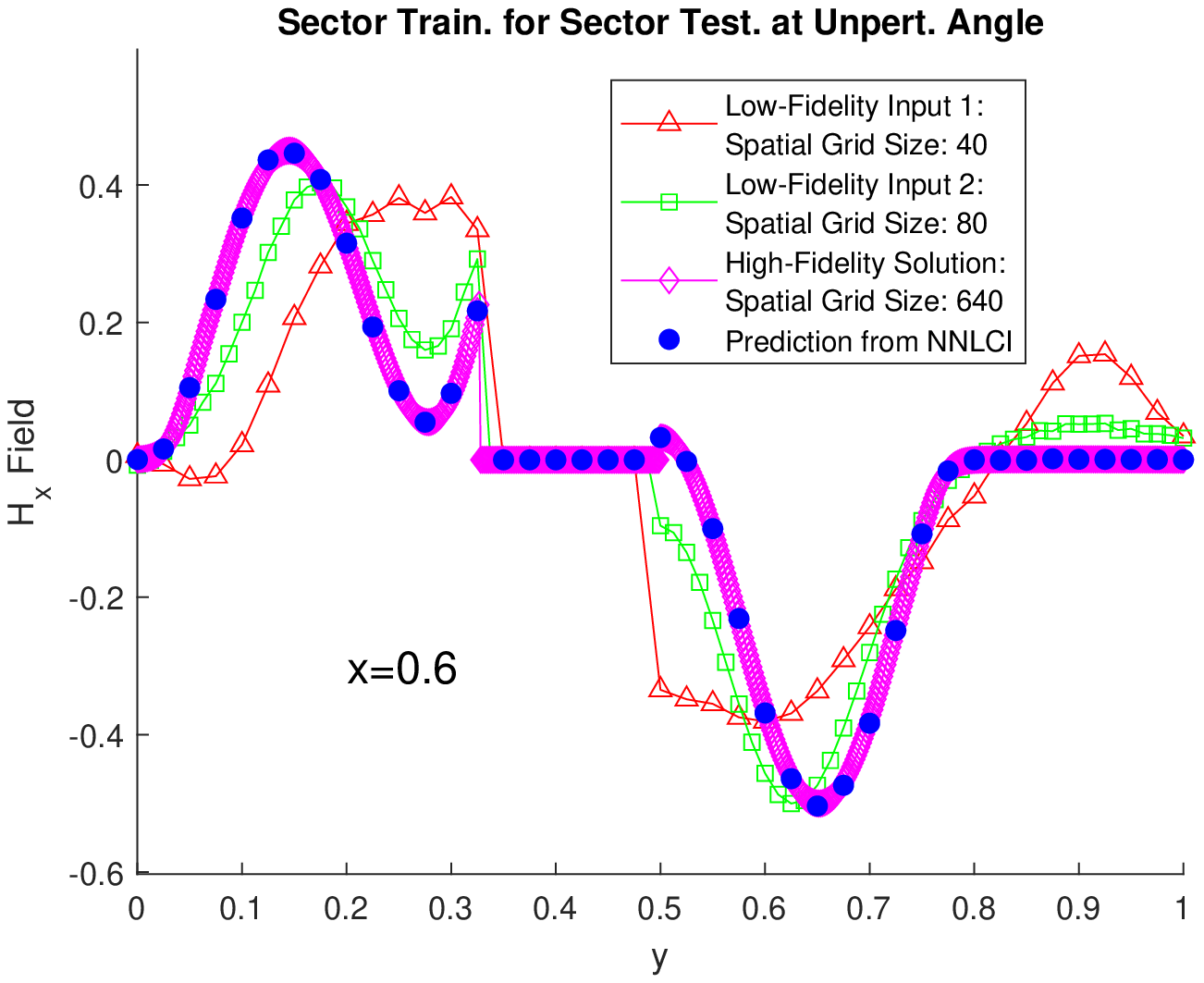}
\captionof{figure}{Cross-sections of the above case at different $y$ values comparing the predicted values with the input and reference values.}\label{fig:90degree_cross}

\begin{tabular}{ |c|c|c|c| }
 \hline
 \multicolumn{4}{|c|}{Sectors from 45 to 135 degrees with 5 degree increments} \\
 \hline
 \multicolumn{2}{|c|}{$\quad$ 10$^\circ$ Increments in Training$\quad$} & \multicolumn{2}{|c|}{15$^\circ$ Increments in Training}\\ 
 \hline
Finer Input Err. & 0.1987 & Finer Input Err. & 0.1991\\
Pred. Err. & 0.0344 & Pred. Err. & 0.0459\\
 \hline
\end{tabular}
\captionof{table}{Accuracy of the predicted solutions to PEC sectors of varying angles}
\label{tbl:sec}
\end{center}

\end{subsubsection}

\begin{subsubsection}{Between circular and square PECs}
We consider the circular PEC as well as the square PEC from subsection 3.1.2. The two PECs are topologically equivalent, yet it is to the best of our knowledge that the corresponding solutions need to be computed separately for each case. To illustrate the effectiveness of our NNLCI-based approach, we generate three different sets of training data for each parameter considered above, namely  all the perturbations (up to $\pm 10 \%$) in frequencies, amplitudes, and radii of the circular PEC. For each training data, we use the unperturbed case for the square PEC as the testing case. We can see from Table \ref{tbl:sfc} that both the frequencies and the radius cases perform similarly well while the amplitude data performs poorly due to the degeneracy of the training data. The results show that our approach of local training {is able to cope with the qualitative differences in the scattering around the corners of the rectangular PEC}.

Conversely, we train on the square PEC subject to all perturbations in the frequencies and side lengths of the rectangle. Then we test on the unperturbed case for the circular PEC, and we observe limited accuracy improvement (Table \ref{tbl:cfs}). One possible explanation is that the  training data for the square PEC cover a small range of the angles at which the incident waves are interacting with the PEC, which may fall short of rendering the predictive power for a broader range of the angles for the circular PEC. On a related note, we leave it for future work to explore the possibility of accurately predicting the cases for arbitrary PEC shapes with only a handful of relatively simple training data. 

\begin{center}

\includegraphics[scale=0.5]{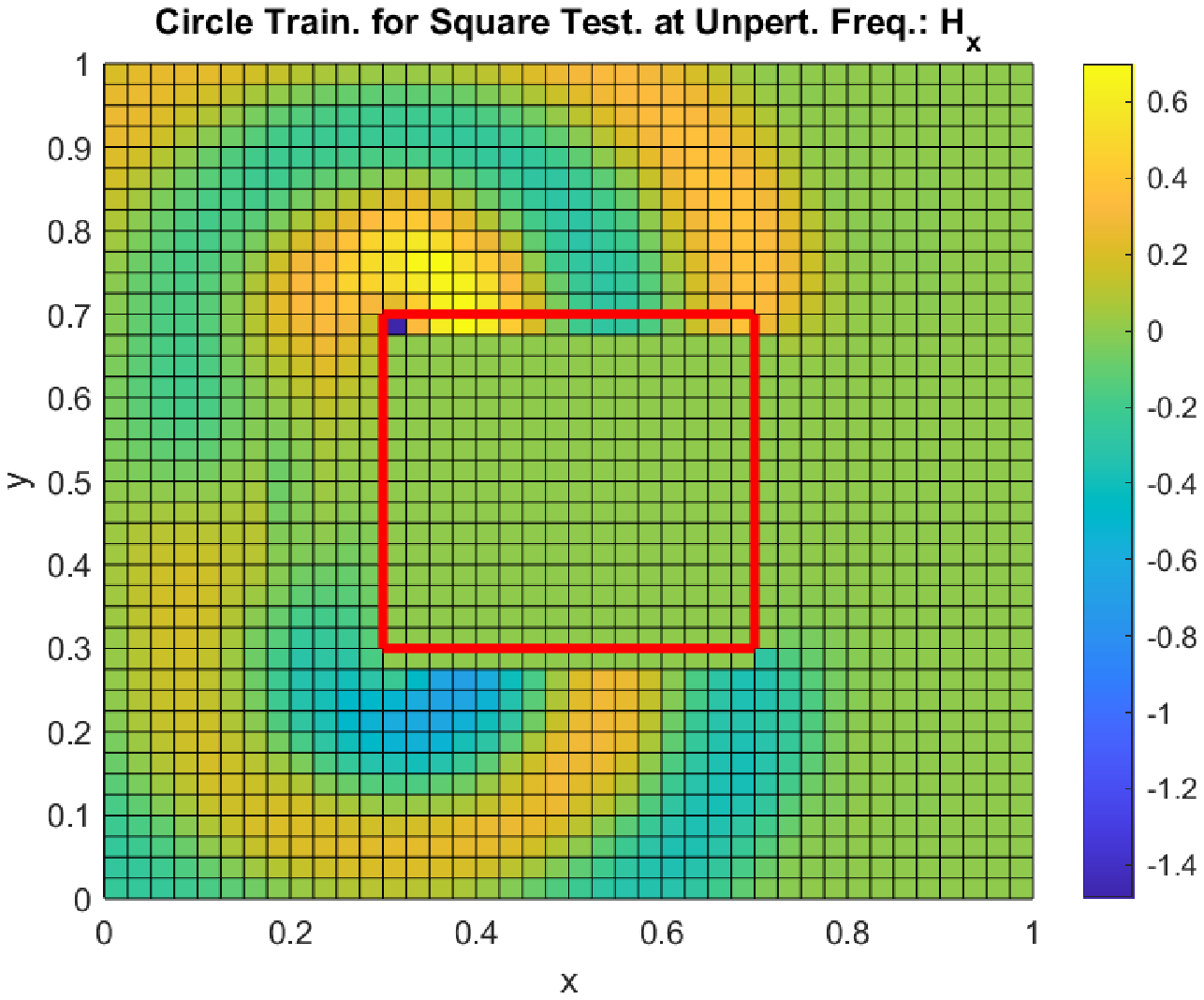}\\
\captionof{figure}{Different shape prediction of $H_x$ scattering off a square in the plane at 0.8 seconds after training the model on all the even frequency perturbations of a circle}
\includegraphics[scale=0.5]{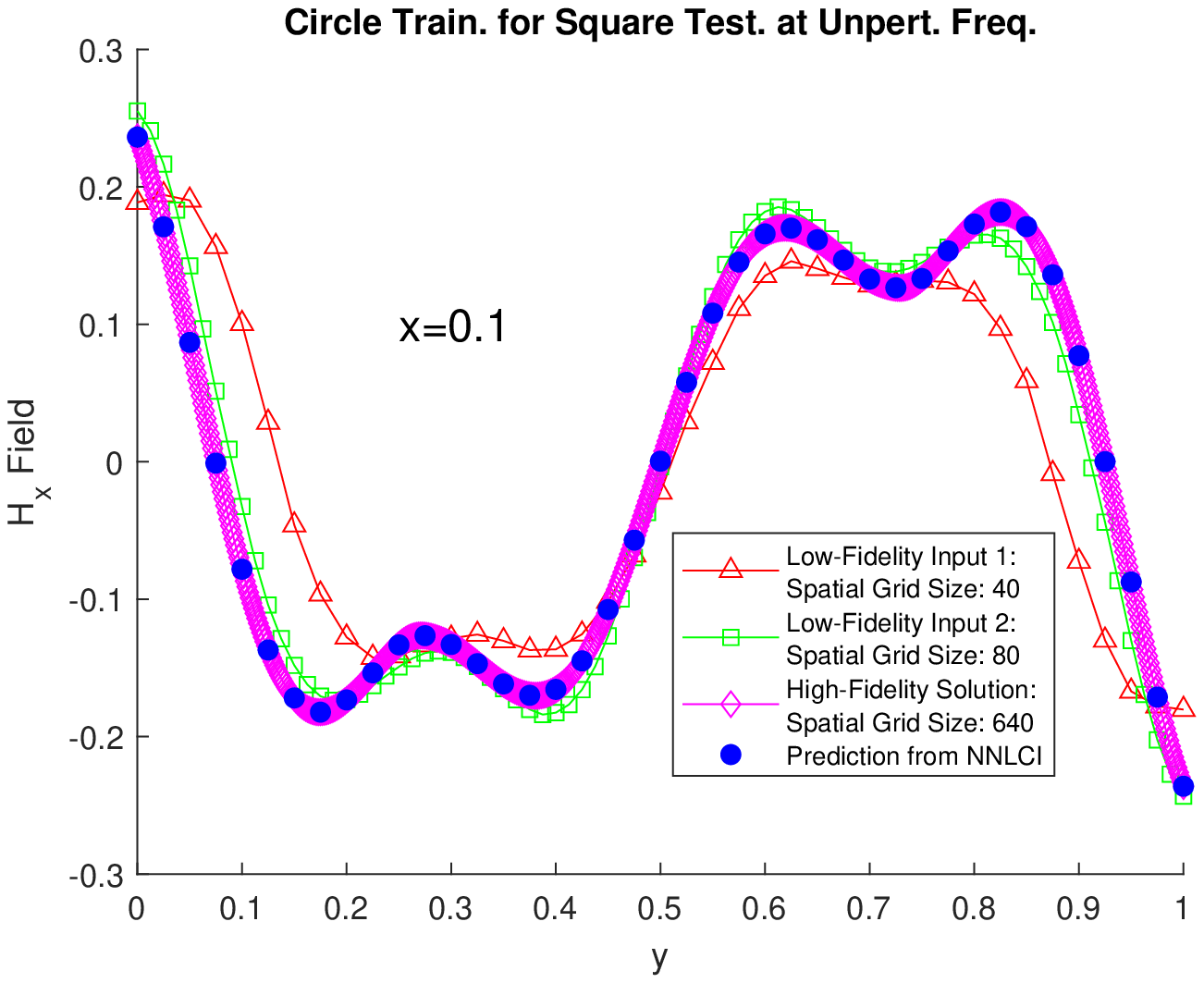}
\includegraphics[scale=0.5]{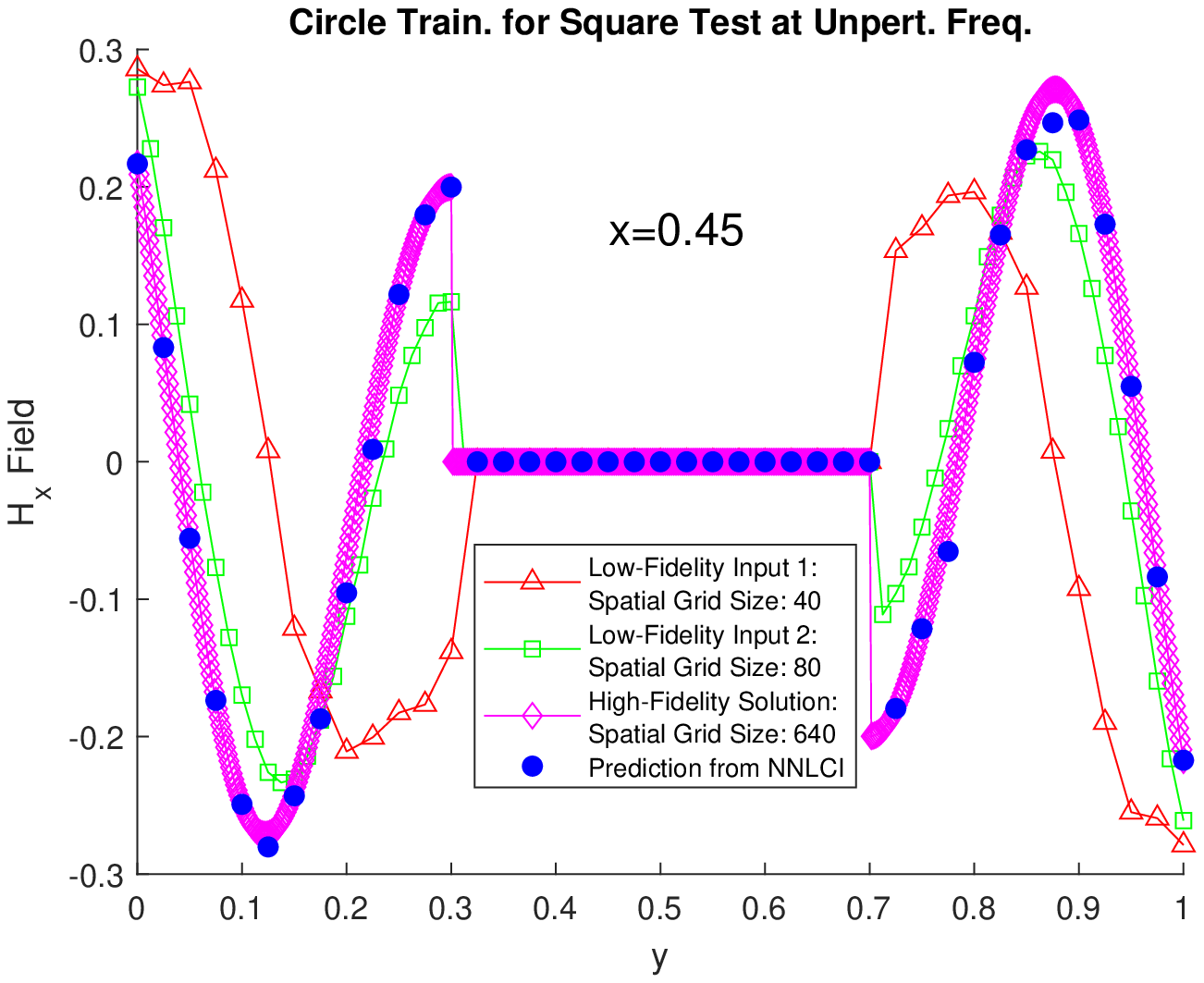}
\captionof{figure}{Cross-sections of the above case at different $y$ values comparing the predicted values with the input and reference values}
\begin{tabular}{ |c|c|c|c|c|c|  }
 \hline
 \multicolumn{6}{|c|}{Predictions on Square PEC from Circular PEC} \\
 \hline
 \multicolumn{2}{|c|}{Circular Frequency} & \multicolumn{2}{|c|}{Circular Amplitude} & \multicolumn{2}{|c|}{Circular Radius}\\
 \hline
Finer Input Err. & 0.1484 & Finer Input Err. & 0.1668 & Finer Input Err. & 0.1480\\
Pred. Err. & 0.0917 & Pred. Err. & 0.2190 & Pred. Err. & 0.0937\\
 \hline
\end{tabular}
\captionof{table}{Accuracy of the predicted solutions to square PEC via training for circular PEC}
\label{tbl:sfc}
\end{center}

\begin{center}

\includegraphics[scale=0.5]{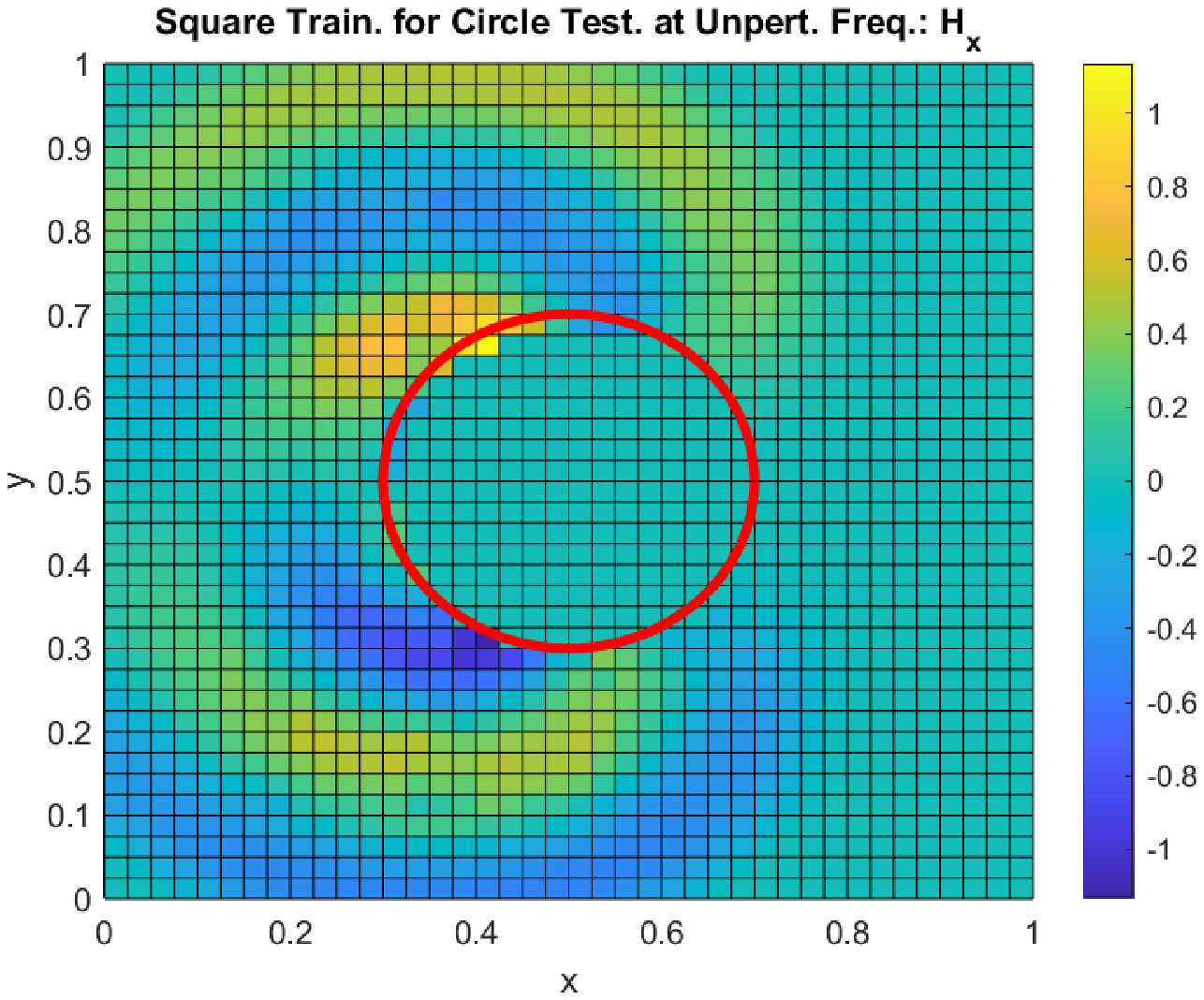}\\
\captionof{figure}{Different shape prediction of $H_x$ scattering off a circle in the plane at 0.8 seconds after training the model on all the even frequency perturbations of a square}
\includegraphics[scale=0.5]{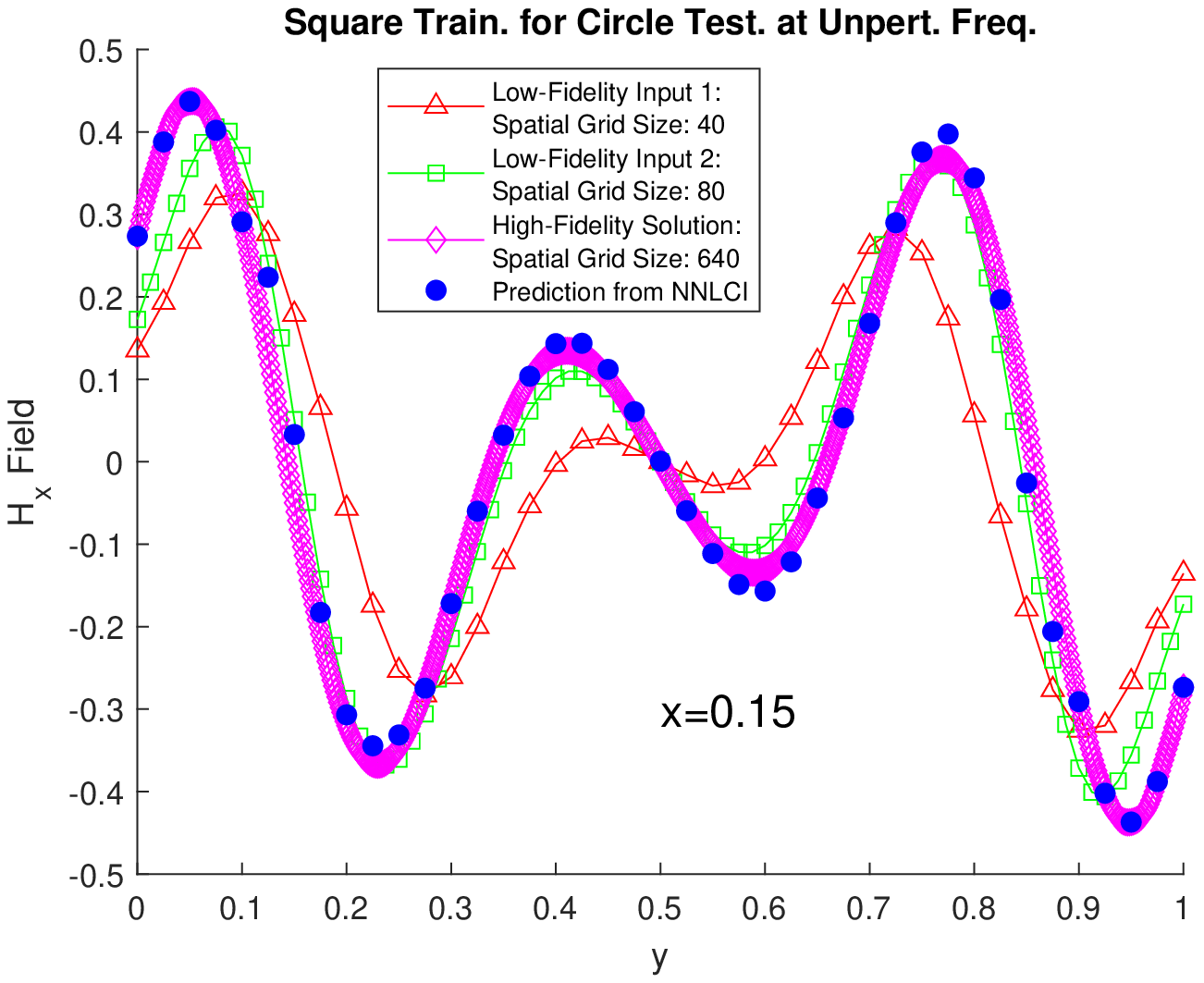}
\includegraphics[scale=0.5]{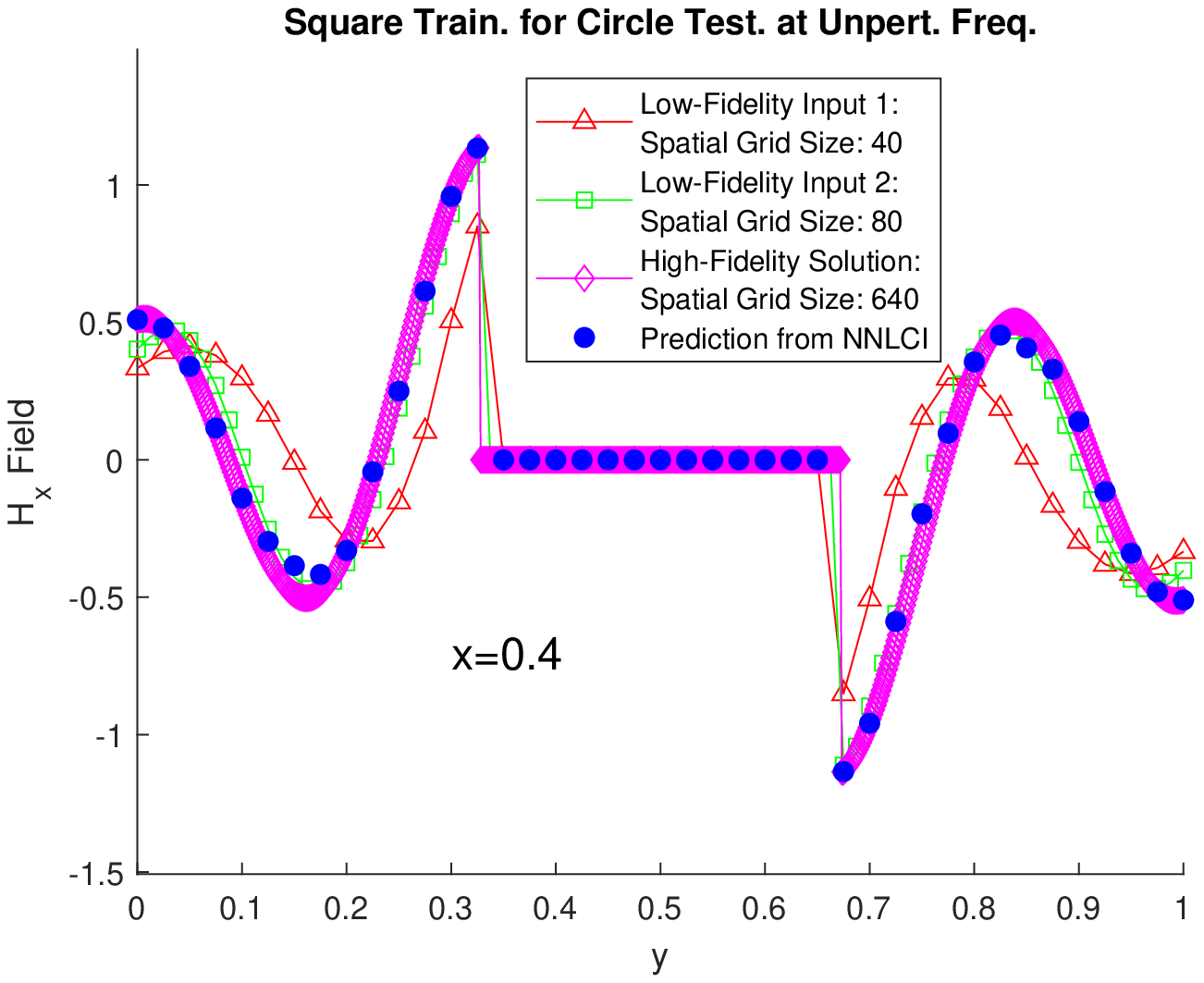}
\captionof{figure}{Cross-sections of the above case at different $y$ values comparing the predicted values with the input and reference values}

\begin{tabular}{ |c|c|c|c|  }
 \hline
 \multicolumn{4}{|c|}{Predictions on Circular PEC from Square PEC} \\
 \hline
 \multicolumn{2}{|c|}{Square Frequency} & \multicolumn{2}{|c|}{Square Side Length} \\
 \hline
Finer Input Err. & 0.1634 & Finer Input Err. & 0.1634\\
Pred. Err. & 0.1442 & Pred. Err. & 0.1274\\
 \hline
\end{tabular}
\captionof{table}{Accuracy of the predicted solutions to circular PEC via training for square PEC}
\label{tbl:cfs}
\end{center}

\end{subsubsection}
\begin{subsubsection}{Quarter Circle Predicts the Square}

In this last test case, we aim at predicting the unperturbed case for the square PEC using the unperturbed (in frequency and amplitude) data for the quarter circle PEC. At first, we considered the training data that consisted of the unperturbed solutions to the original quarter circle as well as to the reflected (about $y = 0.5$) quarter circle, both subject to the incident waves \eqref{eqn:mw} at the terminal time $T = 0.8$. The reflected PEC was included in the training data in order to better capture the corners of the square PEC, yet we observed that our NNLCI was not able to yield accurate predictions.

As a remedy we modify the training data as follows: We shift both the original and reflected quarter circle PECs in the direction of negative $x$-axis by $0.2$ so that the incident waves are in the same phases (with respect to each of the three grid sizes) when they first hit the square and the quarter circle PECs. Then we further enrich the training data by including the three additional training data, namely perturbed frequency cases, 0 and $\pm 2\%$, for the circular PEC. We do this because two training examples (here and in general) are insufficient to give accurate predictions. The result thus obtained is shown in Figs.~\ref{fig:quartertosquare} and \ref{fig:quartertosquare_cross}, and Table \ref{tbl:sqcc}. We remark that the shifted quarter circle PECs may be replaced by the unshifted quarter circles with adjusted terminal times to match the phases of the incident waves for the high-fidelity solution. However, the waves travel at different speeds depending on the grid size due to the dispersive error in our numerical scheme \cite{lee2022ghost} so we must ensure we have training data that covers all three phases. In future work, we wish to carefully examine the propagation of numerical error in our scheme and how that affects the implementation of NNLCI.

\begin{center}

\includegraphics[scale=0.5]{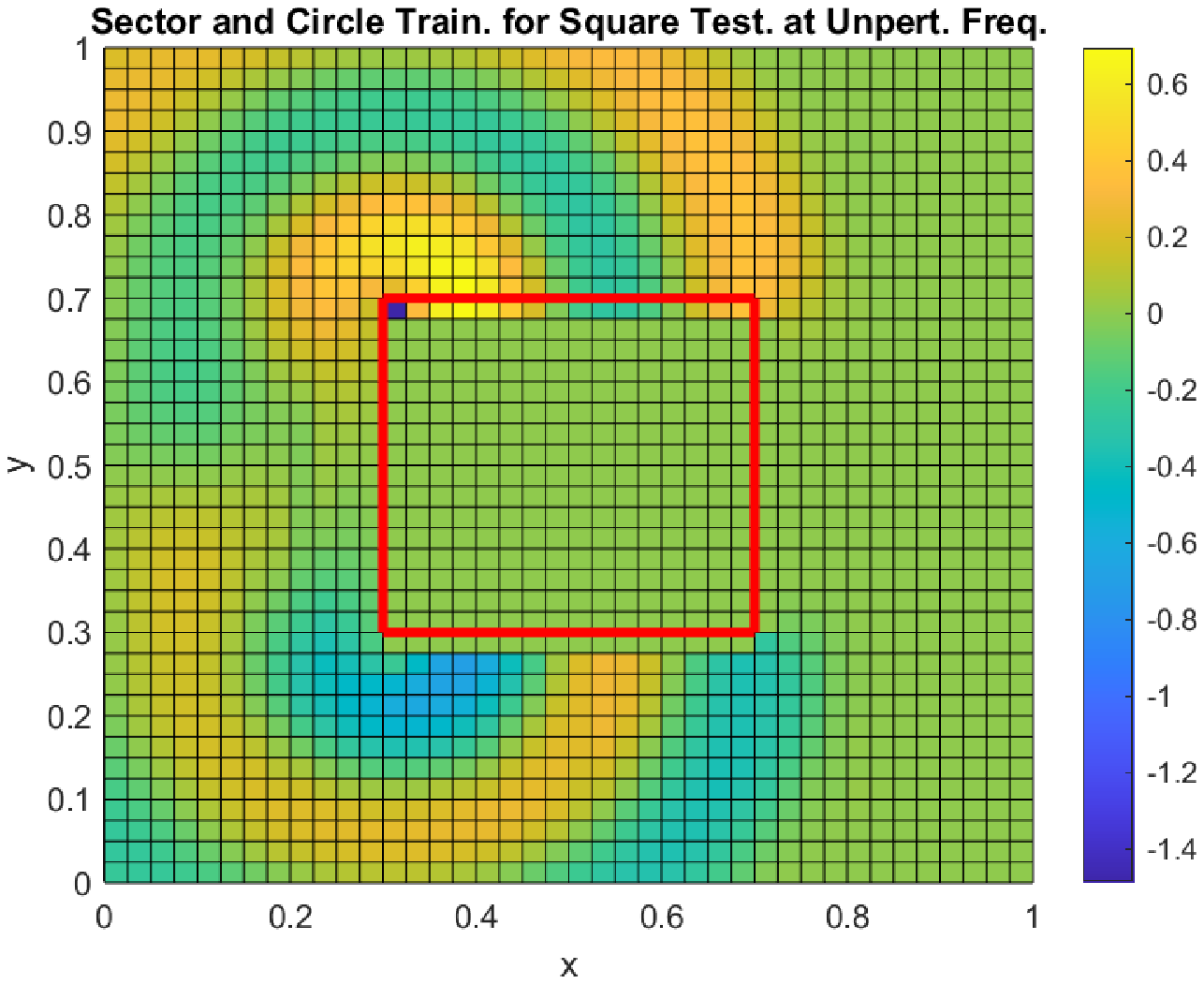}\\
\captionof{figure}{Different shape prediction of $H_x$ scattering off a square in the plane at 0.8 seconds after training the model on a combination of the $90^{\circ}$ sector and circular perturbed frequency data.}\label{fig:quartertosquare}
\includegraphics[scale=0.5]{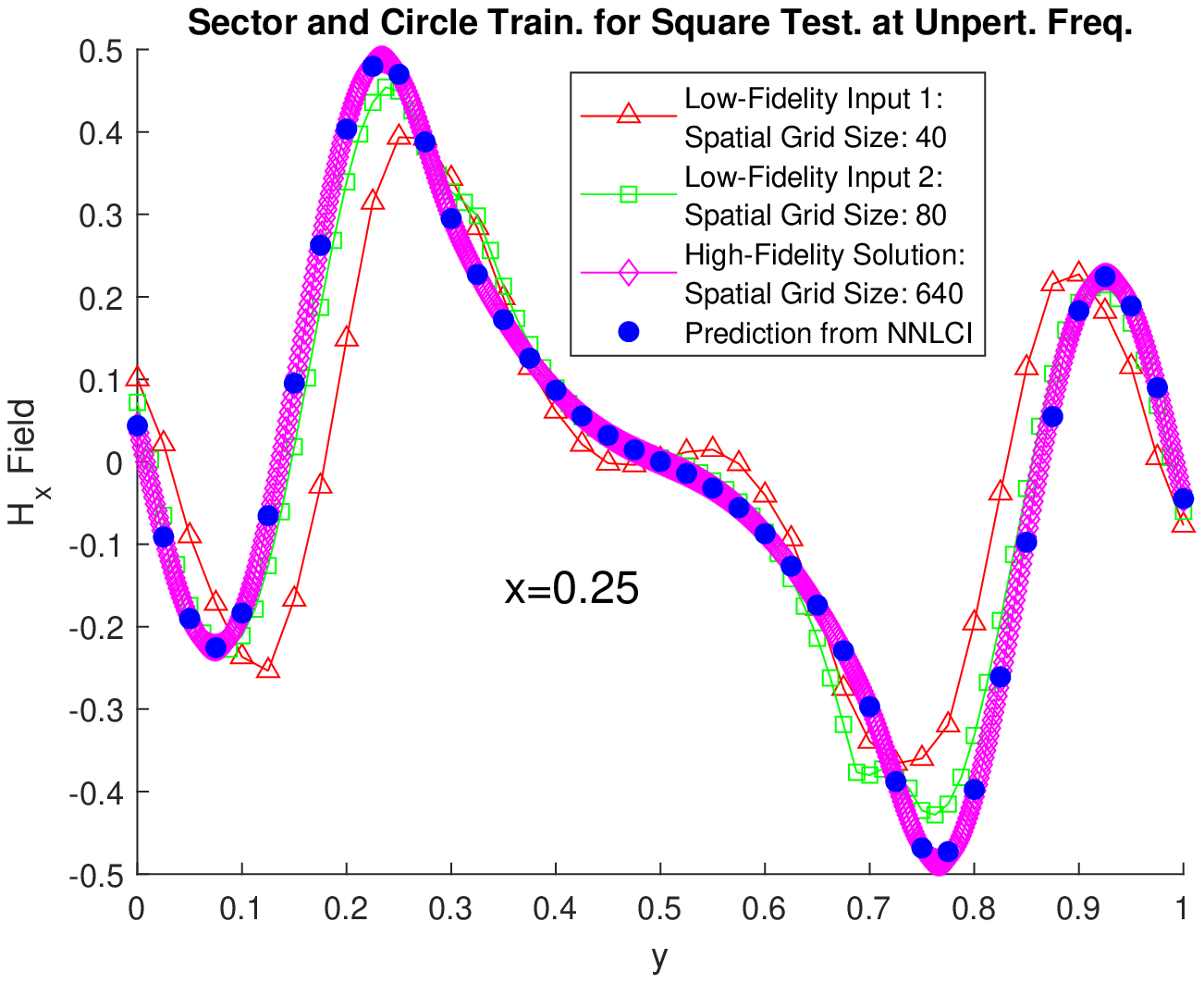}
\includegraphics[scale=0.5]{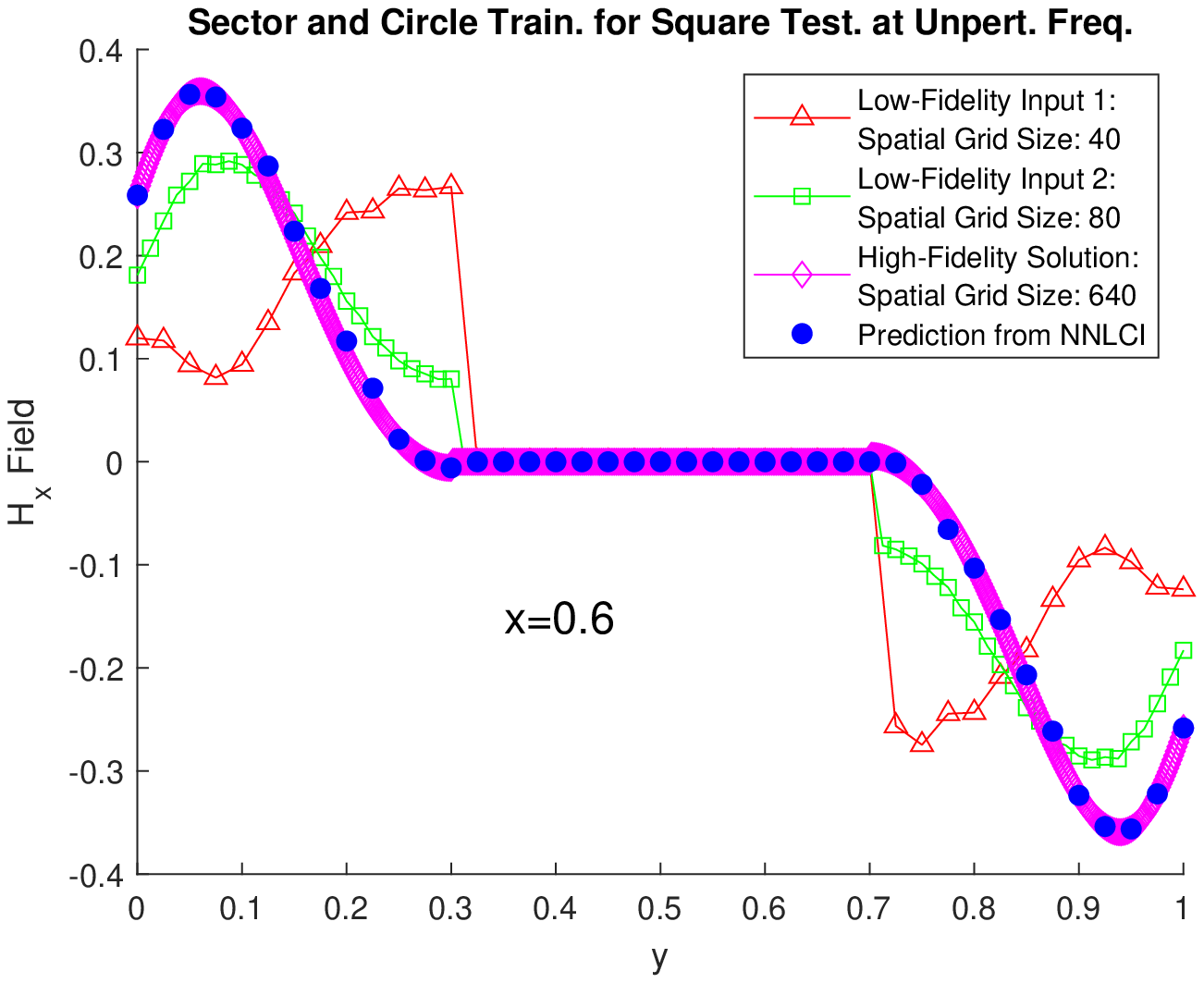}
\captionof{figure}{Cross-sections of the above case at different $y$ values comparing the predicted values with the input and reference values.}\label{fig:quartertosquare_cross}

\begin{tabular}{ |c|c| }
 \hline
 \multicolumn{2}{|c|}{Predicting the Square from the Quarter Circles} \\
 \hline
$\qquad$ Finer Input Err. $\qquad$ & $0.1484$ \\
$\qquad$ Pred. Err. $\qquad$ & $0.0381$ \\
 \hline
\end{tabular}
\captionof{table}{Accuracy of the predicted solutions to square PEC via training on shifted quarter circle PECs and the circular PEC}
\label{tbl:sqcc}
\end{center}

\end{subsubsection}
\end{subsection}

\begin{section}{Conclusion}
\label{sec:4}
In this work, we have applied NNLCI to predict accurate solutions to the scattering of electromagnetic waves around PECs. Once the network is trained, it takes as its inputs two low cost, yet 2nd order accurate, numerical solutions in a converging sequence to generate higher fidelity solutions to PECS of different shapes. Central to the robustness of our approach is the local training of the network which is still effective despite the presence of the PECs. We are interested in theoretical analysis of the NNLCI which may in turn offer insights into choices of suitable grid sizes for the network inputs, as well as possible modifications of the network architecture for enhanced performance. 
\end{section}
\section{Acknowledgements} The authors thank Haoxiang Huang for his help in the coding of NNs.

\bibliographystyle{plain}
\bibliography{Maxwell_NNLCI.bib}

\end{document}